\documentclass[a4paper,12pt]{article}
\begin{document}

\title{Periodic integral operators over Cayley-Dickson algebras and spectra}
\author{Ludkovsky S.V.}
\date{25 March 2012}
\maketitle

\begin{abstract}
Periodic integral operators over Cayley-Dickson algebras related
with integration of PDE are studied. Their continuity and spectra
are investigated. \footnote{key words and phrases: non-commutative
functional analysis, hypercomplex numbers, Cayley-Dickson algebra,
integral operator, spectra, non-commutative integration
 \\
}

\end{abstract}

\section{Introduction}
Integral operators over Cayley-Dickson algebras are useful for
integration of linear and non-linear partial differential equations
\cite{ludvhpde,ludcmft12}.
\par Sections 2 and 3 are devoted to spectra of periodic
integral operators over Cayley-Dickson algebras, that can be used
for analysis of solutions in a bounded domain or of periodic
solutions on the Cayley-Dickson algebra of PDE including that of
non-linear. This is actual, because spectra of operators are used
for solutions of partial differential equations, for example, with
the help of the inverse scattering problem method (see also
\cite{ablsigb}). Moreover, hypercomplex analysis is fast developing
and also in relation with problems of theoretical and mathematical
physics and of partial differential equations
\cite{brdeso,gilmurr,guesprqa}. Cayley-Dickson algebras are used not
only in mathematics, but also in applications
\cite{emch,guetze,girard,krausryan,kravchot}.
\par Analysis over Cayley-Dickson algebras was
developed as well \cite{ludoyst,ludfov,
lujmsalop,lufjmsrf,ludancdnb}. This paper continuous previous
articles and uses their results
\cite{luspraaca,ludspr,ludspr2,ludvhpde,ludcmft12}. Notations and
definitions of papers
\cite{ludoyst,ludfov,lujmsalop,lufjmsrf,ludancdnb} are used below.
The main results of this paper are obtained for the first time.

\section{Periodic integral operators over Cayley-Dickson algebras}

\par {\bf 1. Notation.} Let $X$ be a Banach space over the
Cayley-Dickson algebra ${\cal A}_v$ with $2\le v$, $~v\in \bf N$.
Let also ${\cal A}_w$ be the Cayley-Dickson subalgebra ${\cal
A}_w\subseteq {\cal A}_v$, where $2\le w\le v$. \par The
Cayley-Dickson algebra ${\cal A}_v$ has the real shadow ${\bf
R}^{2^w}$. On this real shadow take the Lebesgue measure $\mu $ so
that $\mu (\prod_{j=1}^{2^w} [a_j,a_j+1])=1$ for each $a_j\in \bf
R$. This measure induces the Lebesgue measure on ${\cal A}_w$
denoted also by $\mu $. A subset $A$ in ${\cal A}_w$ is called $\mu
$-null if a Borel subset $B$ in ${\cal A}_w$ exists so that
$A\subseteq B$ and $\mu (B)=0$. The Lebesgue measure is defined on
the completion ${\cal B}_{\mu }({\cal A}_w)$ of the Borel $\sigma
$-algebra ${\cal B}({\cal A}_w)$ on ${\cal A}_w$ by $\mu $-null
sets.
\par If $1\le p<\infty $, $L^p({\cal A}_w,X)$ will denote the norm
completion of a space of all $\mu $-measurable step functions $f:
{\cal A}_w\to X$ for which the norm
$$(1)\quad \| f \|_p := \sqrt[p]{\int_{{\cal A}_w} \| f(z) \|^p_X \mu (dz)}
$$  is finite, where $$f=\sum_{k=1}^n b_k \chi _{B_k}$$ is a step function so
that $b_k\in X$ and $B_k\in {\cal B}_{\mu }({\cal A}_w)$ for each
$k=1,...,n$; $B_k\cap B_j=\emptyset $ for each $j\ne k$; $ ~ \| x
\|_X$ denotes the norm of a vector $x$ in $X$, certainly, the norm
is non-negative, $0\le \| x \|_X$, $~n$ is a natural number. If
$p=\infty $, the norm is given by the formula:
$$(2)\quad \| f \|_{\infty } := ess_{\mu } ~ \sup_{z\in {\cal A}_w} \| f(z) \|_X
.$$ Then a Banach space $L_q(L^p({\cal A}_w,X))$ of all bounded $\bf
R$-homogeneous ${\cal A}_v$ additive operators $T: L^p({\cal
A}_w,X)\to L^p({\cal A}_w,X)$ is considered. Let $K: {\cal A}_w^2\to
L_q(X)$ be a strongly measurable operator valued mapping, that is a
mapping $g(t,s) := K(t,s)y: {\cal A}_w^2\to X$ is $({\cal B}_{\mu
^2}({\cal A}_w^2), {\cal B}(X))$ measurable for each vector $y\in
X$, i.e. $g^{-1}(Q)\in {\cal B}_{\mu ^2}({\cal A}_w^2)$ for each
Borel subset $Q\in {\cal B}(X)$, where $\mu ^2$ is the Lebesgue
measure on ${\cal A}_w^2$. \par In the paper \cite{ludcmft12} the
following theorem about first order partial differential operators
with variable ${\cal A}_v$ coefficients was demonstrated.
\par {\bf Theorem.} {\it Suppose that a first order
partial differential operator $\Upsilon $ is given by the formula
\par $(i)$ $\Upsilon f= \sum_{j=0}^n  (\partial f/\partial
z_j) {\phi }_j^*(z)$, \\ where ${\phi }_j(z)\ne \{ 0 \}$ for each
$z\in U$ and $\phi _j(z)\in C^0(U,{\cal A}_v)$ for each $j=0,...,n$
such that $Re (\phi _j(z)\phi _k^*(z))=0$ for each $z\in U$ and each
$0\le j\ne k\le n$, where a domain $U$ satisfies Conditions
2.1.1$(D1,D2)$, $ ~ \mbox{}_0z$ is a marked point in $U$, $1<n<
2^v$, $2\le v$. Suppose also that a system $ \{ \phi _0(z),...,\phi
_n(z) \} $ is for $n=2^v-1$, or can be completed by Cayley-Dickson
numbers $\phi _{n+1}(z),...,\phi _{2^v-1}(z)$, such that $(\alpha )$
$ ~ alg_{\bf R} \{ \phi _j(z),\phi _k(z), \phi _l(z) \} $ is
alternative for all $0\le j, k, l \le 2^v-1$ and $(\beta )$
$alg_{\bf R} \{ \phi _0(z),...,\phi _{2^v-1}(z) \} = {\cal A}_v$ for
each $z\in U$. Then a line integral ${\cal I}_{\Upsilon }:
C^0(U,{\cal A}_v)\to C^1(U,{\cal A}_v)$, ${\cal I}_{\Upsilon }f(z)
:= \mbox{}_{\Upsilon }\int_{\mbox{}_0z}^z f(y)dy$ on $C^0(U,{\cal
A}_v)$ exists so that
\par $(ii)$ $\Upsilon {\cal I}_{\Upsilon }f(z) = f(z)$ \\ for each
$z\in U$; this anti-derivative is $\bf R$-linear (or $\bf
H$-left-linear when $v=2$). If there is a second anti-derivative
${\cal I}_{\Upsilon ,2}f(z)$, then ${\cal I}_{\Upsilon }f(z)- {\cal
I}_{\Upsilon ,2}f(z)$ belongs to the kernel $ker (\Upsilon )$ of the
operator $\Upsilon $. }
\par For a first order partial differential operator $\sigma =\Upsilon $ over
${\cal A}_w$ with constant or variable coefficients consider the
antiderivative operator $\mbox{}_{\sigma }\int $ on ${\cal A}_w$.
Put
$$(3)\quad (Bx)(t) := \mbox{}_{\sigma } \int_{-\infty }^{\infty }
K(t,s)x(s)ds $$ whenever this integral converges in the weak sense
as
$$(4)\quad \mbox{}_{\sigma } \int_{- \infty }^{\infty }
u[K(t,s)x(s)yds] := \lim_{a\to a_{\alpha }, b\to b_{\alpha }}
\mbox{}_{\sigma } \int_{\gamma ^{\alpha }|_{[a,b]}} u[K(t,s)x(s)yds]
\in {\cal A}_v$$ for each $y\in X$ and right ${\cal A}_v$ linear
continuous functional $u\in L_r(X,{\cal A}_v)=X_r^*$, where $x\in
L_q(L^p({\cal A}_w,X), L^{p'}({\cal A}_w,X))$ so that
$(xf)(s)=x(s)f(s)$ for each $f\in L^p({\cal A}_w,X)$, $~ x(s)\in
L_q(X)$ for every $s\in {\cal A}_w$, $$\lim_{a\to a_{\alpha }}\gamma
^{\alpha }(t)= \infty \mbox{ and } \lim_{b\to b_{\alpha }}
\gamma^{\alpha }(t) = \infty ,$$ $ a_{\alpha }<b_{\alpha }$, ${\hat
{\cal A}}_w$ is the one-point (Alexandroff) compactification of the
Cayley-Dickson algebra as the topological space, $\infty = {\hat
{\cal A}}_w\setminus {\cal A}_w$, $~\alpha \in \Lambda $. The
integral in Formula $(4)$ reduces to the integral described in \S
4.2.5 \cite{ludancdnb}. Consider a periodic integral kernel
\par $(5)$ $K(t,s)=K(t+p_j\omega _ji_j ,s+p_j\omega _ji_j )$, \\
also $\phi _j(s+p_j\omega _ji_j )=\phi _j(s)$ for $\mu $ almost
every Cayley-Dickson numbers $t, s\in {\cal A}_w$ for all integers
$p_j$, where $\omega _j>0$ is a period by $z_j$,
$~z=z_0i_0+...+z_{2^w-1}i_{2^w-1}\in {\cal A}_w$,
$~z_0,...,z_{2^w-1}\in \bf R$; $~i_0,...,i_{2^w-1}$ denote the
standard generators of ${\cal A}_w$. Suppose that a foliation of
${\cal A}_w$ by paths $\gamma ^{\alpha }$ is so that
\par $(6)$ $K(t,s)=K(\gamma ^{\alpha }(\tau + p^{\alpha }\omega ^{\alpha }),
\gamma ^{\alpha }(\kappa + p^{\alpha }\omega ^{\alpha }))= K(\gamma
^{\alpha }(\tau ),\gamma ^{\alpha }(\kappa ))$, also $\phi _j
(\gamma ^{\alpha }(\kappa + p^{\alpha }\omega ^{\alpha }))= \phi _j
(\gamma ^{\alpha }(\kappa ))$ \\
for $\mu $ almost all $t=\gamma ^{\alpha }(\tau )$ and $s=\gamma
^{\alpha }(\kappa )$ in ${\cal A}_w$ for each $\alpha $ with $\omega
^{\alpha }>0$ and all integers $p^{\alpha }$. Let $S(\omega )$
denote the shift operator
\par $(7)$ $S(\nu )x(t) = x(t+\nu )$ \\ on a Cayley-Dickson number
$\nu \in {\cal A}_w$, where $t\in {\cal A}_w$. \par For example, a
foliation $\{ \gamma ^{\alpha }: ~ \alpha \} $ of ${\cal A}_w$ may
be done by straight lines parallel to ${\bf R}i_0$ indexed by
$\alpha \in \Lambda = \{ z\in {\cal A}_w: ~ Re (z) = 0 \} $.
\par {\bf 2. Definition.} An operator $x\in L_q(L^p({\cal A}_w,X),L^{p'}({\cal A}_w,X))$
so that $(xf)(t)=x(t)f(t)$ and $x(t)\in L_q(X)$ for every $t\in
{\cal A}_w$ we shall call periodic, if
\par $(1)$ $S(\omega _ji_j)x(t)f(t)=x(t)S(\omega _ji_j)f(t)$ \\
for each $t\in {\cal A}_w$ and $f\in L^p({\cal A}_w,X)$ and every
$j$, where $X$ is a Banach space over the Cayley-Dickson algebra
${\cal A}_v$ with $2\le v$. A set $ \{ \omega _j: ~j=0,...,2^w-1 \}
$ will be called a net of periodic values. If for $x$ exists a set
of positive minimal periodic values, then it will be call a set of
periods.

\par {\bf 3. Definition.} Let $Y$ be a Banach space over the
Cayley-Dickson algebra ${\cal A}_v$, where $2\le v$. Put $l_{\infty
}({\bf Z},Y) := \{ x: ~ x: {\bf Z}\to Y, ~ \| x \|_{\infty } :=
\sup_{k\in \bf Z} \| x(k) \|_Y < \infty \} $ and $l_p({\bf Z},Y) :=
\{ x: ~ x: {\bf Z}\to Y, ~ \| x \|_p := [\sum_{k\in \bf Z} \| x(k)
\|_Y^p]^{1/p} < \infty \} $ to be the Banach spaces of norm $\| *
\|_p$ bounded sequences and with values in $Y$, where $\bf Z$
denotes the ring of all integers, $1\le p<\infty $.
\par A sequence $\{ x_n: ~ n\in {\bf N} \} \subset l_{\infty }({\bf
Z},Y)$ is called $c$-convergent to an element $x\in l_{\infty }({\bf
Z},Y)$, if for each integer $k\in \bf Z$ the limit is zero:
$$\lim_{n\to \infty } \| x_n(k)-x(k) \| =0 .$$

\par {\bf 4. Definition.} An operator $B\in L_q(l_{\infty
}({\bf Z},Y))$ will be called $c$-continuous, if an image $\{ Bx_n:
~ n \in {\bf N} \} $ of each $c$-convergent sequence $\{ x_n: ~ n
\in {\bf N} \} $ is a $c$-convergent sequence, where $L_q(X)$ is
written shortly for $L_q(X,X)$, in particular for $X=l_{\infty
}({\bf Z},Y)$, $~Y$ is a Banach space over the Cayley-Dickson
algebra ${\cal A}_v$, $~2\le v$. The family of all $c$-continuous
operators will be denoted by $L_q^c(l_{\infty }({\bf Z},Y))$.

\par {\bf 5. Lemma.} {\it The family $L_q^c(l_{\infty }({\bf
Z},Y))$ from Definition 4 is a closed subalgebra over the
Cayley-Dickson algebra ${\cal A}_v$ in $L_q(l_{\infty }({\bf Z},Y))$
relative to the operator norm topology.}
\par {\bf Proof.} Evidently, the unit operator $I$ is
$c$-continuous. Definitions 3 and 4 imply that $L_q^c(l_{\infty
}({\bf Z},Y))$ is a subalgebra in $L_q(l_{\infty }({\bf Z},Y))$.
Take an arbitrary sequence $B_n$ of $c$-continuous operators
converging to an operator $B\in L_q(l_{\infty }({\bf Z},Y))$
relative to the operator norm. Let $ x_n$ be an arbitrary
$c$-converging sequence to $x$ in $l_{\infty }({\bf Z},Y)$. From the
$c$-continuity of an operator $B_n$ it follows that for each
$\epsilon >0$ and $j\in \bf Z$ there exists a natural number $m$ so
that $\| (B_n(x_k-x) (j) \| <\epsilon /2$ for each $k>m$. On the
other hand, the triangle inequality gives: \par $ \| (B(x_k-x))(j)
\| \le \| ((B_n-B) (x_k-x))(j) \| + \| (B_n(x_k-x))(j) \|$ \par $
\le \| ((B_n-B)( x_k-x))(j) \| +\epsilon /2$. \par The sequence $
B_n$ is norm convergent, hence there exists a natural number $l\in
\bf N$ such that $ \| (B(x_k-x))(j) \| <\epsilon $ for each $k>l$,
since the limit $\lim_{k\to \infty } x_k(j)=x(j)$ exists for each
$j$ with $x\in l_{\infty }({\bf Z},Y)$ and $\sup_k \| (x_k-x)(j)
\|_Y<\infty $. Thus the sequence $Bx_k$ is $c$-convergent: $$\lim_k
Bx_k(j) = Bx(j)\mbox{  for each }j,$$  hence this operator $B$ is
$c$-continuous.

\par {\bf 6. Definitions.} Let $e_k\in l_p({\bf Z},{\cal A}_v)$
be basic elements so that $e_k(j) = \delta _{k,j}$, where $\delta
_{k,j}=0$ for each $j\ne k\in \bf Z$, while $\delta _{j,j}=1$ for
every $j\in \bf Z$. For an operator $B\in L_q(l_p({\bf Z},Y))$ with
$1\le p\le \infty $ and a Banach space $Y$ over the Cayley-Dickson
algebra ${\cal A}_v$ let \par $(1)$ $B_{j,k} y := (Be_ky)(j)\in Y$ \\
for each vector $y\in Y$, where $e_ky\in l_p({\bf Z},Y)$ with
$(e_ky)(j)=\delta _{k,j}y$ for each $j$.
\par This set of operators $ \{ B_{j,k}: ~ j, k\in {\bf Z} \} $ is
called a matrix of an operator $B$.
\par An ${\cal A}_v$ Banach subspace in $l_{\infty }({\bf Z},Y)$ of
all two-sided sequences converging to zero $\lim_{|k|\to \infty }
x(k)=0$ is denoted by $c_0({\bf Z},Y)$.

\par {\bf 7. Lemma.} {\it Suppose that $B\in L_q^c(l_{\infty }({\bf
Z},Y))$ is a $c$-continuous non-zero operator or $B\in L_q(l_p({\bf
Z},Y))$, where $1\le p<\infty $, $Y$ is a Banach space over the
Cayley-Dickson algebra ${\cal A}_v$, $~2\le v$. Then its matrix is
non-zero and bounded. Moreover, $Bc_0({\bf Z},Y)\subset c_0({\bf
Z},Y)$ for each $B\in L_q(l_{\infty }({\bf Z},Y))$.}
\par {\bf Proof.} For each $B\in L_q(l_p({\bf
Z},Y))$ we have the estimate:
\par $ \| B_{j,k} y \| _Y = \| (Be_ky)(j) \| _Y \le \| B e_ky \| _p\le $
\par $\| B \|_{L_q(l_p({\bf Z},Y))} \| e_ky \| _p = \| B
\|_{L_q(l_p({\bf Z},Y))} \| y \| _Y$ \\ for each integers $j, k\in
\bf Z$ and every vector $y\in Y$, hence $ \sup_{j,k} \| B_{j,k} \|
\le \| B \|_{L_q(l_p({\bf Z},Y))}$. Particularly, for $B\in
L_q(l_{\infty }({\bf Z},Y))$ and $x\in c_0({\bf Z},Y)$ this implies
that $$ \| (B(x-\sum_{|k|\le n}e_kx(k))(j) \| \le \| B
\|_{L_q(l_{\infty }({\bf Z},Y))} \| (x-\sum_{|k|\le n}e_kx(k))(j) \|
_Y.$$  But for each $\epsilon >0$ the set $ \{ k: ~ \| x(k) \|_Y
>\epsilon \} $ is finite, since $x\in c_0({\bf Z},Y)$, consequently,
$$\lim_{n\to \infty }\| x-\sum_{|k|\le n}e_kx(k) \| _Y=0$$ and hence
$Bx\in c_0({\bf Z},Y)$.
\par Suppose that $B\in L_q^c(l_{\infty }({\bf Z},Y))$
has a zero matrix $ \{ B_{j,k}: ~ j, k\in {\bf Z} \} $. Take an
arbitrary vector $x\in l_{\infty }({\bf Z},Y)$ and a sequence $\{
x_n: ~ n \in {\bf Z} \} \subset l_{\infty }({\bf Z},Y)\cap c_0({\bf
Z},Y)$ so that it $c$-converges to $x$. In the Banach space
$c_0({\bf Z},Y)$ with norm $ \| x_n \|_{\infty }$ the set of vectors
$\{ e_ky: ~ y \in Y \} $ is everywhere dense. Since $ \{ B_{j,k}: ~
j, k\in {\bf Z} \} =0$, the restriction $B|_{(c_0({\bf Z},Y))}$ is
zero and $Bx_n=0$ for each $n$. Thus the sequence $Bx_n$ does not
converge to $Bx$. This produces the contradiction. Therefore, the
matrix $ \{ B_{j,k}: ~ j, k\in {\bf Z} \} $ is non-zero. \par In the
space $l_p({\bf Z},Y)$ with $1\le p<\infty $ a subset of finite
two-sided sequences in $Y$ is dense, hence a non-zero operator $B\in
L_q(l_p({\bf Z},Y))$ has a non-zero matrix $ \{ B_{j,k}: ~ j, k\in
{\bf Z} \} $.

\par {\bf 8. Definition.} Let $\Lambda $ be a Banach algebra with unit
over the Cayley-Dickson algebra ${\cal A}_v$, $~2\le v$. A
subalgebra $\Xi $ over ${\cal A}_v$ with unit is called saturated if
each element of $\Xi $ invertible in $\Lambda $ is invertible in
$\Xi $ as well.

\par {\bf 9. Remark.} Henceforth in this section it will be
denoted for short by $l({\bf Z},Y)$ each of the Banach spaces
$c_0({\bf Z},Y)$ and $~l_p({\bf Z},Y)$ with $p\in [1,\infty ]$ if
something other will not specified. Henceforward operators from
$L_q^c(l_{\infty }({\bf Z},Y))$ will be considered.
\par Let $g\in l_{\infty }({\bf Z},{\cal A}_v)$, then an operator
$D(g)$ is defined by the formula:
\par $(1)$ $(D(g)x)(k) := g(k)x(k)$ \\
for each $x\in l_{\infty }({\bf Z},Y)$. The family of all such
operators is denoted by $\cal E$. The family of all operators $A\in
L_q^c(l_{\infty }({\bf Z},Y))$ quasi-commuting with each $D(g)\in
\cal E$ is denoted by $L_q^s(l_{\infty }({\bf Z},Y))$, that is
$$(2)\quad \mbox{}^jA \mbox{ }^kD(g) =(-1)^{\kappa (j,k)} \mbox{ }^kD(g)
\mbox{ }^jA$$ for each $j, k =0,1,2,...$ (see \S \S II.2.1
\cite{ludopalglamb}, 2.5 and 2.23 \cite{ludunbnormopla12}).
\par Evidently, this subalgebra $L_q^s(l_{\infty }({\bf Z},Y))$
is saturated in the operator algebra $L_q^c(l_{\infty }({\bf Z},Y))$
with unit operator $I$ as the unit of these algebras.

\par {\bf 10. Definition.} Take a Cayley-Dickson number of absolute value
one $M\in {\bf S}_v := \{ z: ~ z\in {\cal A}_v, ~ |z|=1 \} $ and put
${\cal M} = \{ M^n: ~ n \in {\bf Z} \} \in l_{\infty }({\bf Z},{\cal
A}_v)$, where $2\le v$.
\par An operator $B\in L_q^c(l_{\infty }({\bf Z},Y))$ is called
diagonal if it quasi-commutes with $D({\cal M})\in \cal E$ for each
$M\in {\bf S}_v$, where $Y$ is a Banach space over the
Cayley-Dickson algebra ${\cal A}_v$, $~2\le v$. The family of all
diagonal operators we denote by $L_q^d(l_{\infty }({\bf Z},Y))$.

\par {\bf 11. Proposition.} {\it The family of all diagonal operators
$L_q^d(l_{\infty }({\bf Z},Y))$ is a saturated subalgebra in the
algebra $L_q^c(l_{\infty }({\bf Z},Y))$.}
\par {\bf Proof.} If $A, B\in L_q^d(l_{\infty }({\bf Z},Y))$, then
equations 9$(2)$ for $A$ and $B$ and $D({\cal M})$ imply
$$(1)\quad (\mbox{}^jA+\mbox{}^jB) \mbox{ }^kD({\cal M})= \mbox{}^jA\mbox{ }^kD({\cal M})
+\mbox{}^jB \mbox{ }^kD({\cal M})=(-1)^{\kappa (j,k)} \mbox{
}^kD({\cal M}) (\mbox{ }^jA+\mbox{}^jB)$$ for each $j, k =0,1,2,...$
and every $M\in {\bf S}_v$ due to distributivity of the operator
multiplication (see also \S 3.3.1). If $z\in {\cal A}_v$, then $zA$
has components:
$$(2)\quad \mbox{}^j(zA)=\sum_{s,l; ~ i_si_l=i_j} (z_si_s\mbox{ }^lA +
(-1)^{\kappa (s,l)} z_li_l \mbox{ }^sA)$$ for each $j$, where as
usually $z=z_0i_0+z_1i_1+...$ with $z_0, z_1,...\in \bf R$.
Therefore, $zA$ quasi-commutes with each $D({\cal M})$ for each
$z\in {\cal A}_v$. Analogously it can be demonstrated for $Az$. Thus
$L_q^d(l_{\infty }({\bf Z},Y))$ is an algebra.
\par An inverse $A^{-1}\in L_q^d(l_{\infty }({\bf Z},Y))$
of an operator $A\in L_q^d(l_{\infty }({\bf Z},Y))$ quasi-commuting
with $D({\cal M})$ also quasi-commutes with $D({\cal M})$, since
$D({\cal M})$ is invertible and $(AD({\cal M}))^{-1} = (D({\cal
M}))^{-1}A^{-1} = D({\cal M}^*)A^{-1}$ and $(D({\cal M})A)^{-1} =
A^{-1} D({\cal M}^*)$, $~Ax=y$ implies $x=A^{-1}y$, $$ ~ I =
(AD({\cal M}))^{-1} (AD({\cal M})) = $$ $$(AD({\cal M}))^{-1}
\sum_{j; ~ i_si_l=i_j} [\mbox{}^sA\mbox{ }^lD({\cal M}) +
(-1)^{\kappa (s,l)} \mbox{ }^lA\mbox{ }^sD({\cal M})]$$ $$=
(AD({\cal M}))^{-1} \sum_{j; ~ i_si_l=i_j} [ \mbox{}^sD({\cal
M})\mbox{ }^lA + (-1)^{\kappa (s,l)} \mbox{ }^lD({\cal M})\mbox{
}^sA ]$$ $$ = \sum_{j; ~ i_pi_q=i_j^*; ~ i_si_l=i_j} \{
[\mbox{}^pD({\cal M}^*)\mbox{ }^qA^{-1} +(-1)^{\kappa (p,q)}\mbox{
}^qD({\cal M}^*)\mbox{ }^pA^{-1}]$$ $$[ \mbox{}^sD({\cal M})\mbox{
}^lA + (-1)^{\kappa (s,l)} \mbox{ }^lD({\cal M})\mbox{ }^sA ] \} ,$$
where ${\cal M}^*$ corresponds to $M^* = \tilde M$.

\par {\bf 12. Definition.} A sequence $g\in l_{\infty }({\bf Z},{\cal A}_v)$
is called periodic of period $k\in \bf N$ if $g(n+k)=g(n)$ for each
integer $n\in \bf Z$, where $Y$ is a Banach space over the
Cayley-Dickson algebra ${\cal A}_v$, $~2\le v$.
\par {\bf 13. Lemma.} {\it Diagonal operators from $L_q^d(l_{\infty
}({\bf Z},Y))$ quasi-commute with operators of multiplication on
periodic sequences $g\in l_{\infty }({\bf Z},{\cal A}_v)$, where
$2\le v$.}
\par {\bf Proof.} Let $A\in L_q^d(l_{\infty }({\bf Z},Y))$ be a
diagonal operator and let $g\in l_{\infty }({\bf Z},{\cal A}_v)$ be
a periodic sequence of period $k\in \bf N$. Put $\theta _j := \exp
(2\pi {\bf i}j/k)$, where $~j= 0, 1,..., k-1$, $ ~ {\bf i}$ is an
additional purely imaginary generator so that ${\bf i}^2=-1$, ${\bf
i}i_l=i_l{\bf i}$ for each $l\ge 0$. \par A minimal real algebra
with basis of generators $i_0, i_1,..., i_{2^v-1}, {\bf i}, {\bf
i}i_1,..., {\bf i}i_{2^v-1}$ and their relations as above is the
complexification $({\cal A}_v)_{{\bf C}_{\bf i}}$ of the
Cayley-Dickson algebra ${\cal A}_v$, where ${\bf C}_{\bf i} = {\bf
R}\oplus {\bf R}{\bf i}$. Then $g$ can be presented in the form:
$$(1)\quad g(n) = \sum_{j=0}^{k-1} c_j \theta _j^n, \mbox{
 where}$$
$$(2)\quad c_j= \frac{1}{k} \sum_{n=0}^{k-1} g(n) {\tilde \theta
}_j^n$$ are $({\cal A}_v)_{{\bf C}_{\bf i}}$ Fourier coefficients
for each $j=0,1,...,k-1$. Indeed, ${\bf i}({\bf i}x_ji_j)=-x_ji_j$
for each $j$ and $x_j\in Y_j$, while $(a+b)c=ac+bc$ and
$c(a+b)=ca+cb$ for each $c\in ({\cal A}_v)_{{\bf C}_{\bf i}}$ and
$a, b \in Y\oplus Y{\bf i}$, where we put ${\bf i}x=x{\bf i}$ for
each $x\in Y$. Put $ \| x+y{\bf i} \| ^2 = \| x \| ^2 + \| y \| ^2$
for each $x, y \in Y$ and $x+y{\bf i} \in Y\oplus Y{\bf i}$.
Therefore,
$$\sum_{j=0}^{k-1}c_j \theta _j^m = \frac{1}{k}\sum_{j=0}^{k-1}
[\sum_{n=0}^{k-1} (\sum_l g_l(n)i_l) {\tilde \theta }_j^n] \theta
_j^m$$
$$=\sum_l g_l(n) \frac{1}{k} \sum_{j=0}^{k-1}
\sum_{n=0}^{k-1} i_l \theta _j^{m-n}= g(m),$$ since the
multiplication in the Cayley-Dickson algebra is distributive, where
$g(m)=\sum_l g_l(m)i_l$ with $g_l(m)\in \bf R$ for each
$l=0,1,2,...$, $$\frac{1}{k} \sum_{j=0}^{k-1} \exp( 2\pi j{\bf
i}(m-n)/k) =\delta _{n,m}.$$ Therefore, the diagonal operator has
the decomposition
$$(D(g)x)(n) = g(n)x(n) =
(\sum_{j=0}^{k-1} (c_jI)D({\hat \theta }_j)x)(n)$$ that is
$$(3)\quad D(g)=\sum_{j=0}^{k-1} (c_jI)D({\hat \theta }_j)$$ again due to
distributivity of the multiplication in the Cayley-Dickson algebra
${\cal A}_v$ (see \S 2.5 \cite{ludunbnormopla12}). On  the other
hand, $A$ and $D({\hat \theta }_j)$ quasi-commute:
$$\mbox{}^jA\mbox{ }^kD({\hat \theta }_j)= (-1)^{\kappa (j,k)}
 \mbox{ }^kD({\hat \theta }_j)\mbox{ }^jA$$ for each $j, k$
and from Proposition 11 it follows that $A$ and $D(g)$
quasi-commute, where $\mbox{}^kD({\hat \theta }_j)(n) \in {\bf
C}_{\bf i}i_k$ for each $n\in \bf Z$.
\par {\bf 14. Corollary.} {\it Diagonal operators from
$L_q^d(l_{\infty }({\bf Z},Y))$ quasi-commute with multiplication
operators $D(g)$ on periodic sequences $g\in l_{\infty }({\bf
Z},{\cal A}_v)$, where $2\le v$.}
\par {\bf 15. Definition.} If $x\in l_p({\bf Z},Y)$, where
$p\in [1,\infty ]$, $Y$ is a Banach space over the Cayley-Dickson
algebra ${\cal A}_v$, $~2\le v$, then its support is $supp ~ x := \{
n: ~n\in {\bf Z}; ~ x(n)\ne 0 \} $.
\par {\bf 16. Lemma.} {\it Let $B\in L_q^d(l_{\infty }({\bf Z},Y))$
be a diagonal operator and let $x\in l_{\infty }({\bf Z},Y)$ be with
finite support, where $Y$ is a Banach space over the Cayley-Dickson
algebra ${\cal A}_v$, $~2\le v$. Then $supp ~Bx \subseteq supp ~x.$}
\par {\bf Proof.} If the support of $x$ is finite, then there exists
a natural number $N\in \bf N$ so that $supp ~x \subseteq [-N,N]$.
Consider natural numbers $m>N$ and $n\in [-N,N]$. Put $g(n+2mk) = 1$
if $n\in supp ~x$, while $g(n+2mk) = 0$ for $n\in [-N,N]\setminus
supp ~x$, where $k\in \bf Z$, hence $g\in l_{\infty }({\bf Z},{\cal
A}_v)$ is of period $2m$. Therefore, $D(g)x=x$ by the construction
of $g$. In view of Corollary 14 $G(g)$ and $B$ quasi-commute. On the
other hand, $g(l)=0$ for each $l\in [-m,m]\setminus supp ~x$,
consequently, $(Bx)(l)=0$, since $Bx=BD(g)x$. Thus $(Bx)(l)=0$ for
each $l\notin supp ~x$, since $m>N$ is arbitrary.
\par {\bf 17. Theorem.} {\it The algebras $L_q^s(l_{\infty }({\bf
Z},Y))$ and $L_q^d(l_{\infty }({\bf Z},Y))$ over the Cayley-Dickson
algebra ${\cal A}_v$ coincide, where $2\le v$.}
\par {\bf Proof.} From the definitions it follows that the algebra
$L_q^s(l_{\infty }({\bf Z},Y))$ is a subalgebra of the algebra
$L_q^d(l_{\infty }({\bf Z},Y))$ of diagonal continuous operators.
Therefore, it remains to prove the inclusion $L_q^d(l_{\infty }({\bf
Z},Y))\subseteq L_q^s(l_{\infty }({\bf Z},Y)).$ Consider arbitrary
continuous diagonal operator $B\in L_q^d(l_{\infty }({\bf Z},Y))$
and an element $g\in l_{\infty }({\bf Z},{\cal A}_v)$. Take $x\in
l_{\infty }({\bf Z},Y)$ so that $x_n=0$ for each $|n|>N$, where $N$
is a natural number. We have $\mbox{}^jD(g)\mbox{ }^kBx =
(-1)^{\kappa (j,k)} \mbox{ }^kB\mbox{}^jD(g)x$ for each $j, k$.
Extend the sequence $g$ periodically to $h$ so that $h(-N)=g(N+1)$
and $h(m)=g(m)$ for each $m\in [-N,N]$, consequently,
$(D(g)x)(m)=(D(h)x)(m)$ for each $m\in [-N,N]$. In view of Corollary
14 $$(\mbox{}^jB\mbox{ }^kD(g)x)(m)= (\mbox{}^jB\mbox{ }^kD(h)x)(m)
=(-1)^{\kappa (j,k)}(\mbox{ }^kD(h)\mbox{}^jBx)(m)$$ $$=
(-1)^{\kappa (j,k)}(\mbox{ }^kD(g)\mbox{}^jBx)(m)$$ for each $j, k$
and $m\in [-N,N]$. Applying Lemma 16 for each $|m|>N$ we get
$(Bx)(m) =0$ and hence $$(\mbox{}^jB\mbox{ }^kD(g)x)(m)=(-1)^{\kappa
(j,k)}(\mbox{ }^kD(g)\mbox{}^jBx)(m)$$ for each $j, k=0,1,2,...$ and
every integer number $m\in \bf N$. Thus for sequences with finite
supports this theorem is accomplished. But a set of all sequences
with finite support is dense in $c_0({\bf Z},Y)$ and in $l_p({\bf
Z},Y)$ for every $1\le p<\infty $. Therefore, the statement of this
theorem is valid on these spaces. In accordance with Lemma 7 the
conjecture spreads on a $c$-continuous operator $B$ from $c_0({\bf
Z},Y)$ on the entire Banach space $l_{\infty }({\bf Z},Y)$.
\par {\bf 18. Definition.} Let $g= (G^m: ~ m\in {\bf Z})$ be a sequence
belonging to the Banach space $l_{\infty }({\bf Z},L_q(Y))$, where
$Y$ is a Banach space over the Cayley-Dickson algebra ${\cal A}_v$,
$~2\le v$. An operator $D(g)\in L_q^c(l_{\infty }({\bf Z},Y))$ will
be defined by the formula: $(D(g)x)(m) = G^m x(m)$ for each integer
number $m$. A set of all (left) multiplication operators on bounded
operator valued sequences forms an algebra over the Cayley-Dickson
algebra ${\cal A}_v$, which will be denoted by $L^b_q(l_{\infty
}({\bf Z},Y)))$. \par An operator $B\in L_q^c(l_{\infty }({\bf
Z},Y))$ is called a $(k,n)$ ribbon operator with $k\in \bf N$ and
$n\in \bf Z$ if $B_{s,l} =0_Y$ is the zero operator from $L_q(Y)$
for each $|s-l+n|\ge  k$, where $s, l\in \bf Z$. Their family is
denoted by $L_q^{(k,n)}(l_{\infty }({\bf Z},Y))$. A $(k,0)$ ribbon
operator is called $k$-ribbon (single ribbon for $k=1$).
\par Two propositions follow immediately from the latter definition.
\par {\bf 19. Proposition.} {\it An operator $B\in L_q^c(l_{\infty }
({\bf Z},Y))$ is $(1,0)$ ribbon if and only if it is a
multiplication operator on operator valued sequence.}
\par {\bf 20. Proposition.} {\it The algebra $L_q^{(k,0)}(l_{\infty }({\bf Z},Y))$
is the saturated subalgebra of the algebra $L_q^c(l_{\infty }({\bf
Z},Y))$ over the Cayley-Dickson algebra ${\cal A}_v.$}
\par {\bf 21. Theorem.} {\it Let $B\in L_q^c(l_{\infty }({\bf Z},Y))$,
where $Y$ is a Banach space over the Cayley-Dickson algebra ${\cal
A}_v$, $~2\le v$. An operator $B$ is $(1,0)$ ribbon if and only if
$B$ is diagonal.}
\par {\bf Proof.} Suppose that $B$ is a diagonal operator. Take a
vector $x\in Y$ and an element $y^s=e_sx\in l_{\infty }({\bf Z},Y),$
where $y^s(k)=\delta _{s,k}x$, hence $B_{m,s}(By)(m)\in  Y$, where
$s, m \in \bf Z$, $B_{s,m}$ are elements of the matrix of $B$. For
$M\in {\bf S}_v$ matrix elements of the operator $D({\cal
M}^*)BD({\cal M})$ are prescribed by the formula: \par $(D({\cal
M}^*)BD({\cal M})y^s)(m)=M^{-m} (BM^sy^s)(m)$ \\ for each $s, m \in
{\bf Z}$. \par Take any purely imaginary generator $i_p$ of the
Cayley-Dickson algebra and put $M=i_p$ with $p\ge 1$. As an operator
$B$ is diagonal, the equalities follow: $$(\mbox{}^kBy^s)(m) =
(-1)^{\kappa (p,k)\eta (s)} i_p^{s-m} (\mbox{}^kBy^s)(m)$$ for each
$m, s$, where $\eta (s)=0$ for $s$ even, while $\eta (s)=1$ for $s$
odd. This implies that $B_{m,s}=0_Y$ for $s\ne m$ and
$B_{m,s}=B_{m,m}$ for $m=s$. Thus the operator $B$ is $(1,0)$
ribbon.
\par The inverse conjecture follows from Lemma 7 and Definition 18.
\par {\bf 22. Corollary.} {\it An operator $B\in L_q^c(l_{\infty
}({\bf Z},Y))$ is diagonal if and only if $B$ is an operator of
multiplication on a bounded operator valued sequence in $L_q(Y)$.}

\par {\bf 23. Definition.} An operator $B\in L_q(l_{\infty }({\bf
Z},Y)))$ is called uniformly $c$-continuous, if a mapping $\breve{B}
: {\bf S}^1\to L_q(l_{\infty }({\bf Z},Y\oplus Y{\bf i})))$ is
continuous relative to the operator norm topology on $L_q(l_{\infty
}({\bf Z},Y\oplus Y{\bf i})))$ and a topology on ${\bf S}^1 := \{ z:
~ z\in {\bf C}_{\bf i}; ~ |z| =1 \} $ induced by the norm on the
complex field ${\bf C}_{\bf i}$, where $2\le v$, $~ \breve{B}(M) :=
D({\cal M})BD({\cal M}^*)$ for each $M\in {\bf S}^1$, $ ~ D({\cal
M})$ is a diagonal (left) multiplication operator on a sequence
${\cal M}(k) = M^k$, $~Y$ is a Banach space over ${\cal A}_v$. The
family of all uniformly $c$-continuous operators is denoted by
$L_q^{uc}(l_{\infty }({\bf Z},Y))$ and is supplied with the uniform
operator norm topology
$$\| B \|_u := \sup_{M\in {\bf S}^1} \| \breve{B}(M) \| .$$
The real field is the center of the Cayley-Dickson algebra ${\cal
A}_v$ with $v\ge 2$, hence the generator ${\bf i}$ can be realized
as the real matrix ${{~0 ~ 1}\choose{-1~0}}$. If $X$ and $Y$ are two
${\cal A}_v$ vector spaces and $B: X\to Y$ is a real homogeneous
${\cal A}_v$ additive operator, then it has a natural extension
${\bf B}: X_{\bf i}\to Y_{\bf i}$ so that ${\bf B}(a+b{\bf i})
=(Ba)+(Bb){\bf i}$ for each vectors $a, b\in X$, where $X_{\bf i}$
is obtained form $X$ by extending the algebra ${\cal A}_v$ to
$({\cal A}_v)_{{\bf C}_{\bf i}}$ and $X_{\bf i}$ can be presented as
the direct sum $X_{\bf i} = X\oplus X{\bf i}$ of two ${\cal A}_v$
vector spaces $X$ and $X{\bf i}$. It is convenient to denote ${\bf
B}$ also by $B$ on $X_{\bf i}$.
\par {\bf 24. Proposition.} {\it The family $L_q^{uc}(l_{\infty }({\bf
Z},Y))$ is a normed algebra over the Cayley-Dickson algebra ${\cal
A}_v$ (see \S 23). If an operator $B\in L_q^c(l_{\infty }({\bf
Z},Y))$ is $k$-ribbon, then this operator $B$ is uniformly
$c$-continuous.}

\par {\bf Proof.} In algebra $L_q(X)$ for a Banach space $X$ over
the Cayley-Dickson algebra ${\cal A}_v$ the operator norm satisfies
the inequality: $ \| AB \| \le \| A \| \| B \| $, particularly for
$X=l_{\infty }({\bf Z},Y)$ or $X=l_{\infty }({\bf Z},Y\oplus Y{\bf
i})$. Therefore, for each uniformly $c$-continuous operators $A$ and
$B$ the inequality follows: $$\| AB \|_u =\sup_{M\in {\bf S}^1} \|
D({\cal M}) AB D({\cal M}^*) \| = \sup_{M\in {\bf S}^1} \| D({\cal
M}) AD({\cal M}^*)D({\cal M})B D({\cal M}^*) \| $$ $$\le \sup_{M\in
{\bf S}^1} \| D({\cal M}) AD({\cal M}^*) \| \sup_{M\in {\bf S}^1} \|
D({\cal M}) BD({\cal M}^*) \| = \| A \|_u \| B \|_u,$$ since set
theoretic composition of operators is associative and $D({\cal M})
D({\cal M}^*) = I$ for each $M\in {\bf S}^1$ (see \cite{schafb}). On
the other hand, $\| B \| _u\ge \| B \| $ for each $B\in
L_q^{uc}(l_{\infty }({\bf Z},Y))$, since $D({\cal I})=I$ and $1\in
{\bf S}^1$, where ${\cal I}=(...,1,1,...)$ corresponds to $1$. \par
For a $k$ ribbon operator $B$ a finite sequence $\mbox{}_nB$ of
single ribbon operators exists with $|n|\le k$ so that $$B=
\sum_{n=-k}^k \mbox{ }_nB S(n),$$ where $S(n)$ denotes a shift
operator on $n$, $~ (S(n)g)(m) := g(m+n)$ for each $n, m \in \bf Z$
and $g\in l_{\infty }({\bf Z},Y)$. Therefore, the equality follows:
$$(\breve{B}(M) x)(j) = \sum_{n=-k}^k [(M^{j+n}I) \mbox{
}_nB ((M^*)^nI) x](j+n) ,$$ but $\| [(M^{j+n}I) \mbox{ }_nB
((M^*)^nI) y \| \le \| \mbox{ }_nB y \| $ for each $y\in Y$ and
$M\in {\bf S}^1$, consequently, $ \| B \|_u \le \sum_{n=-k}^k \|
\mbox{ }_nB  \|  <\infty $.

\par {\bf 25. Lemma.} {\it Let $B\in L_q^{uc}(l_{\infty }({\bf Z},Y))$ be a uniformly
$c$ continuous operator and $(B_{s,p})$ be its matrix, where $s, p
\in \bf Z$, where $Y$ is a Banach space over the Cayley-Dickson
algebra ${\cal A}_v$, $~2\le v$. Then for each $M\in {\bf S}^1$
matrix elements of an operator $\breve{B}(M)$ have the form:
$$\breve{B}(M)_{s,p} = (M^sI)B_{s,p}(M^{-p}I) .$$}
\par {\bf Proof.} In accordance with Definition 6 the equalities are valid:
\par $\breve{B}(M)_{s,p} = (D({\cal M}) B D({\cal M}^*))_{s,p}x =
((D({\cal M}) B D({\cal M}^*)) e_px)(s)$\par $ =
((M^sI)B(M^{-p}I)e_px)(s)= ((M^sI)B_{s,p}(M^{-p}I)e_px)(s)$ \\
for all integer numbers $s, p\in \bf Z$ and for each vector $x\in
Y$, since $M\in {\bf S}^1$ implies $|M|^2=MM^*=M^*M=1$ and hence
$M^{-1}=M^*$.

\par {\bf 26. Proposition.} {\it Let $A, B\in L_q^{uc}(l_{\infty }({\bf Z},Y))$
be two uniformly $c$-continuous operators, where $Y$ is a Banach
space over the Cayley-Dickson algebra ${\cal A}_v$, $~2\le v$. Then
$\breve{(AB)}(M) = \breve{A}(M)\breve{B}(M)$ for each $M\in {\bf
S}^1$.}
\par {\bf Proof.} The algebra of operators relative to the
set-theoretic composition is evidently associative (see also
\cite{schafb}), hence \par $\breve{(AB)}(M) = D({\cal M}) AB D({\cal
M}^*) =D({\cal M}) A D({\cal M}^*)D({\cal M}) B D({\cal M}^*)$\par $
= \breve{A}(M)\breve{B}(M)$ for each $M\in {\bf S}^1$, \\ since
$MM^*=M^*M=1$ and hence $D({\cal M}) D({\cal M}^*) = D({\cal M}^*)
D({\cal M}) $.

\par {\bf 27. Proposition.} {\it Let $B\in L_q^{uc}(l_{\infty }({\bf Z},Y))$
be a uniformly $c$-continuous invertible operator, where $Y$ is a
Banach space over the Cayley-Dickson algebra ${\cal A}_v$, $~2\le
v$. Then an operator $\breve{B}(M)$ is invertible so that
$\breve{B^{-1}}(M) = (\breve{B})^{-1}(M)$ for each $M\in {\bf
S}^1$.}
\par {\bf Proof.} Applying Proposition 26 with $A=B^{-1}$ one gets
$(\breve{B})^{-1}(M) = (D({\cal M}) B D({\cal M}^*))^{-1} = D({\cal
M}^*)^{-1} B^{-1} D({\cal M})^{-1} = D({\cal M}) B^{-1} D({\cal
M}^*) =\breve{B^{-1}}(M)$.

\par {\bf 28. Lemma.} {\it Let $B\in L_q^{uc}(l_{\infty }({\bf Z},Y))$
be a uniformly $c$-continuous operator, where $Y$ is a Banach space
over the Cayley-Dickson algebra ${\cal A}_v$, $~2\le v$. Then there
exists an equivalent norm on $Y$ relative to an initial one so that
$ \| B \| \ge \| \breve{B}(M) \| $ for each $M\in {\bf S}^1$.}
\par {\bf Proof.} The norm on the Cayley-Dickson algebra ${\cal
A}_v$ satisfies the inequality $ | ab | \le | a | | b | $ for each
$a, b \in {\cal A}_v$ with $2\le v$. Particularly, for $v\le 3$ the
norm on ${\cal A}_v$ is multiplicative. \par Two norms $\| * \| $
and $ \| * \| '$ on a Banach space $Y$ are called equivalent if two
positive constants $0<c_1\le c_2<\infty $ exist so that $ c_1 \| x
\| \le \| x \| ' \le c_2 \| x \| $ for each vector $x\in Y$. Then $
\| ax \| \le | a | \| x \| $ for each $a\in {\cal A}_v$ and $x\in Y$
up to a topological isomorphism of Banach spaces, i.e. up to an
equivalence of norms on $Y$, since $\| tx_ji_j \| = |t| \| x_j \| =
\| t x_j \| $ for each $x_j\in Y_j$ and $j=0, 1, 2,...$. Indeed, the
multiplication of vectors on numbers ${\cal A}_v\times Y\ni
(a,x)\mapsto ax\in Y$ is continuous relative to norms on ${\cal
A}_v$ and $Y$. Therefore, $ \| \breve{B}(M) x \| = \| D({\cal M}) B
D({\cal M}^*) x \| \le \| D({\cal M}) \| \| B \| \| D({\cal M}^*)\|
\| x \| $, consequently, $ \| B \| \ge \| \breve{B}(M) \| $ for each
$M\in {\bf S}^1$, since $ |M|^n=|M^n|=1$ for each integer $n$ and
hence $ \| D({\cal M}) \| =1$.

\par {\bf 29. Definition.} Let $C_s({\bf S}^1, L_q^{uc}(l_{\infty
}({\bf Z},Y\oplus Y{\bf i}))$ denote a Banach space of continuous
bounded mappings from ${\bf S}^1$ into $L_q^{uc}(l_{\infty }({\bf
Z},Y\oplus Y{\bf i}))$, where $Y$ is a Banach space over the
Cayley-Dickson algebra ${\cal A}_v$, $~2\le v$.
\par {\bf 30. Corollary.} {\it There exists an equivalent norm on a Banach
space $Y$ over the Cayley-Dickson algbera ${\cal A}_v$ such that the
mapping $F: L_q^{uc}(l_{\infty }({\bf Z},Y)) \to C_s({\bf S}^1,
L_q^{uc}(l_{\infty }({\bf Z},Y\oplus Y{\bf i}))$ given by the
formula $F(B)(M)=\breve{B}(M)$ for each $M\in {\bf S}^1$ is $\bf R$
linear and ${\cal A}_v$ additive and isometric operator.}
\par {\bf Proof.} This follows by combining Proposition 24 and
Lemma 28.

\par {\bf 31. Lemma.} {\it Let $B\in L_q^{uc}(l_{\infty }({\bf Z},Y))$. Then
an operator valued mapping $\breve{B}: {\bf S}^1\to
L_q^{uc}(l_{\infty }({\bf Z},Y))$ has the Fourier series of the
form:
$$(1)\quad \breve{B} \sim \sum_{n= - \infty }^{\infty } M^n\mbox{ }_n\breve{B}
$$ for each $M\in {\bf S}^1$, where
$$(2)\quad \mbox{ }_n\breve{B} = \frac{1}{2\pi } \int_0^{2\pi }
e^{-nt{\bf i}} \breve{B}(e^{t{\bf i}}) dt$$ are Fourier
coefficients. Moreover, each operator $\mbox{ }_n\breve{B}S(-n)$ is
diagonal.}
\par {\bf Proof.} From the definition of the uniformly $c$-continuous
operator it follows that the restriction of  $\breve{B}$ on ${\bf
S}^1$ is continuous, since a mapping $\breve{B}: {\bf S}^1\to
L_q^{uc}(l_{\infty }({\bf Z},Y)$ is continuous. The function
$e^{t{\bf i}}$ has the period $2\pi $. Therefore, integrals $(2)$
exist for every $n$ and a formal Fourier series $(1)$ can be
written. On the other hand, the algebra $alg_{\bf R} (i_k,i_l)$ is
associative for every $k$ and $l=0,1,2,...$. Therefore, using the
distributivity law in the algebra $({\cal A}_v)_{{\bf C}_{\bf i}}$
we deduce that $$M^{\mp k} B_{k,l} M^{\pm l} \sum_p x_pi_p= M^{\mp
k} B_{k,l} \sum_p M^{\pm l} x_pi_p,$$  where $M\in {\bf S}^1$,
$x_p\in Y_p$ for each $ p =0, 1, 2,...$. The algbera $alg_{\bf R}
({\bf i},i_p,i_k)$ is associative for each $p, k$, consequently, the
inversion formula $(1)$ is valid, since
$$\frac{1}{2\pi } \int_0^{2\pi } e^{-nt{\bf i}} e^{mt{\bf i}} \mbox{
}_m\breve{B} x dt = \delta_{m,n} \mbox{ }_m\breve{B} x$$ for each
vector $x\in Y$. Then one gets \par $D({\cal M}) \mbox{ }_n\breve{B}
S(-n) D({\cal M}^*) = D({\cal M}) \mbox{ }_n\breve{B}D({\cal M}^*)
 D({\cal M}) S(-n) D({\cal M}^*)$\par $ = (M^nI)\mbox{
 }_n\breve{B}(M^{-n}I) (M^nI)S(-n)(M^{-n}I) = (M^nI)\mbox{
 }_n\breve{B}S(-n)(M^{-n}I)$ \\ for each $M\in {\bf S}^1$
 and every integer $n$, since the product of diagonal operators is
 diagonal.

 \par {\bf 32. Theorem.} {\it The algebra $L_q^{uc}(l_{\infty }({\bf
 Z},Y))$ is the saturated subalgebra in $L_q^c(l_{\infty }({\bf
 Z},Y))$, where $Y$ is a Banach space over the Cayley-Dickson
algebra ${\cal A}_v$, $~2\le v$.}
\par {\bf Proof.} The algebraic $\bf R$ linear ${\cal A}_v$ additive
embedding $L_q^{uc}(l_{\infty }({\bf
 Z},Y))\hookrightarrow L_q^c(l_{\infty }({\bf Z},Y))$ follows from
 the defintions. If a uniformly $c$-continuous operator $B$ on $l_{\infty }({\bf
 Z},Y)$ is invertible in the algebra $L_q^c(l_{\infty }({\bf Z},Y))$, then
 the mapping $\breve{B}^{-1}$ is
 continuous from ${\bf S}^1$ into $L_q^c(l_{\infty }({\bf Z},Y))$
due to Proposition 27. Thus $B^{-1}\in L_q^{uc}(l_{\infty }({\bf
 Z},Y))$.

\par {\bf 33. Definition.} An operator $B\in L_q^c(l_{\infty }({\bf
Z },Y))$ is called periodic of period $n$ on the Banach space
$l_{\infty }({\bf Z },Y)$ over the Cayley-Dickson algebra ${\cal
A}_v$ with $2\le v$ if $S(n)B=BS(n)$, where $n$ is a natural number,
$(S(n)x)(k) = x(n+k)$ for each vector $x\in l_{\infty }({\bf Z},Y)$
and every integer $k$. A set of $n$ periodic operators will be
denoted by $L_q^{n,per}(l_{\infty }({\bf Z },Y))$.

\par {\bf 34. Proposition.} {\it A set $L_q^{n,per}(l_{\infty }({\bf Z },Y))$
of $n$ periodic operators on $l_{\infty }({\bf Z },Y)$ is a closed
saturated subalgebra in the algebra $L_q^c(l_{\infty }({\bf Z
},Y))$, where $Y$ is a Banach space over the Cayley-Dickson algebra
${\cal A}_v$, $~2\le v$.}
\par {\bf Proof.} If operators $A, B\in L_q^c(l_{\infty }({\bf Z
},Y))$ commute with $S(n)$, then \par $[(\alpha I)A + B(\beta
I)]S(n)= (\alpha I)S(n)A + B S(n) (\beta I)$\par $ = S(n) (\alpha
I)A + S(n)B(\beta I)=S(n)[(\alpha I)A + B(\beta I)]$ \\ for each
Cayley-Dickson numbers $\alpha , \beta \in {\cal A}_v$ and
\par $ABS(n)=AS(n)B=S(n)AB$. \par Thus $L_q^{n,per}(l_{\infty }({\bf Z
},Y))$ is an algebra over the Cayley-Dickson algebra ${\cal A}_v$.
From Definition 33 the algebraic $\bf R$ linear ${\cal A}_v$
additive embedding $L_q^{n,per}(l_{\infty }({\bf Z
},Y))\hookrightarrow L_q^c(l_{\infty }({\bf Z },Y))$ follows. \par
The relation $S(n)B-BS(n)=0$ defines a closed subset in
$L_q^c(l_{\infty }({\bf Z },Y))$, since $S(n)$ is the bounded
continuous operator on $l_{\infty }({\bf Z},Y)$ and the mapping
$f(B) := S(n)B-BS(n)$ is continuous from $L_q^c(l_{\infty }({\bf Z
},Y))$ into itself $L_q^c(l_{\infty }({\bf Z },Y))$. Thus
$L_q^{n,per}(l_{\infty }({\bf Z },Y))$ is the closed subalgebra.
\par If an $n$ periodic operator $B\in L_q^{n,per}(l_{\infty }({\bf Z },Y))$
is invertible in $L_q^c(l_{\infty }({\bf Z },Y))$, then
$B^{-1}S(n)=S(n)B^{-1}$, since $BS(n)=S(n)B$ $\Leftrightarrow $
$S(n)= B^{-1}S(n)B$ $\Leftrightarrow $ $S(n)B^{-1} = B^{-1}S(n)$.
Thus $B^{-1}\in L_q^{n,per}(l_{\infty }({\bf Z },Y))$ and hence the
subalgebra $L_q^{n,per}(l_{\infty }({\bf Z },Y))$ is saturated in
$L_q^c(l_{\infty }({\bf Z },Y))$.

\par {\bf 35. Lemma.} {\it An operator $B\in L_q^c(l_{\infty }({\bf Z },Y))$ is
$n$ periodic if and only if its matrix satisfies the condition
$B_{k+n,l+n} = B_{k,l}$ for each integers $k$ and $l$, where $Y$ is
a Banach space over the Cayley-Dickson algebra ${\cal A}_v$, $~2\le
v$.}
\par {\bf Proof.} Suppose that $B\in L_q^{n,per}(l_{\infty }({\bf Z },Y))$
and $(B_{k,l})$ is its matrix. Then $B_{k,l}x = (Be_lx)(k) =
(S(-n)BS(n) e_lx)(k) = (S(-n) B e_{l+n}x)(k) = (Be_{l+n}x)(k+n) =
B_{k+n,l+n}x$ for each vector $x\in Y$ and integers $k$ and $l$.
\par Vise versa if $B\in L_q^c(l_{\infty }({\bf Z },Y))$ is a
$c$-continuous operator, then it has a matrix $(B_{k,l})$ by Lemma
7. Then the condition $B_{k+n,l+n} = B_{k,l}$ for each integers $k$
and $l$ implies $B_{k+n,l+n}x= (Be_{l+n}x)(k+n)= (S(-n) B
e_{l+n}x)(k) = (S(-n)BS(n) e_lx)(k)= (Be_lx)(k)$, consequently,
$S(-n)BS(n)=B$ and hence $BS(n)=S(n)B$. Thus the operator $B$ is
$n$-periodic.

\par {\bf 36. Corollary.} {\it Suppose that $B$ is a $c$-continuous
operator $B\in L_q^c(l_{\infty }({\bf Z },Y))$. Then $B$ is
$n$-periodic and diagonal if and only if it is a (left)
multiplication operator on a stationary $n$-periodic sequence in
$L_q(Y)$.}
\par {\bf Proof.} This follows from Theorem 17 and Lemma 35.

\par {\bf 37. Definition.} A function ${\hat B}: {\bf S}^1 \to
L_q(Y\oplus Y{\bf i})$ prescribed by the formula ${\hat B} (M)x :=
\bigoplus_{j=0}^{n-1} B(D({\cal M})x)(j) $ will be called the
Fourier transform of an $n$-periodic operator $B\in
L_q^{n,per}(l_{\infty }({\bf Z },Y))$, where $Y$ is a Banach space
over the Cayley-Dickson algebra ${\cal A}_v$, $~2\le v$. We put $$
\| {\hat B} (M)x \| := \max_{j=0}^{n-1} \| B(D({\cal M})x)(j) \| .$$
\par By $C_s({\bf S}^1,(L_q(Y\oplus Y{\bf i}))^n)$
will be denoted the Banach space of all bounded continuous mappings
$G: {\bf S}^1\to L_q(Y\oplus Y{\bf i})$ supplied with the norm $$ \|
G \| := \sup_{M\in {\bf S}^1} \max_{j=0}^{n-1} \| \mbox{}_jG(D({\cal
M})) \| ,$$ where $~\| A \| $ denotes a norm of an operator $A\in
L_q(Y)$, $~Y$ is a Banach space over the Cayley-Dickson algebra
${\cal A}_v$, $~2\le v$, $~G = \bigoplus_{j=0}^{n-1} \mbox{}_jG$
with $\mbox{}_jG\in L_q(Y\oplus Y{\bf i})$ for every $j$.

\par {\bf 38. Lemma.} {\it Let $B\in L_q^{n,per}(l_{\infty }
({\bf Z },Y))$ be a periodic operator, where $Y$ is a Banach space
over the Cayley-Dickson algebra ${\cal A}_v$, $~2\le v$. Then ${\hat
B} \in C_s({\bf S}^1,(L_q(Y\oplus Y{\bf i}))^n)$.}
\par {\bf Proof.} A uniform space $C_s({\bf S}^1,L_q(Y\oplus Y{\bf i}))$
is complete for a Banach space $Y$ over the Cayley-Dickson algebra
${\cal A}_v$. Take an arbitrary vector $x\in Y$ and a complex number
$K\in {\bf S}^1$ and a sequence $\mbox{}_kM\in {\bf S}^1$ converging
to $K$. The Banach spaces $(L_q(Y\oplus Y{\bf i}))^n$ and
$\bigoplus_{j=0}^{n-1}L_q(Y\oplus Y{\bf i})$ are isometrically
isomorphic when supplied with the corresponding norms, since $n$ is
a natural number, where \par $(1)$ $\| A \| = \sup_{0\le j\le n-1}
\| \mbox{}_jA \| $ \\ for each $A=(\mbox{}_0A,...,\mbox{}_{n-1}A)\in
\bigoplus _{j=0}^{n-1} L_q(Y\oplus Y{\bf i})$, also \par $(2)$ $ \|
x \| = \sup_{0\le j\le n-1} \| \mbox{}_jx \| $ \\ for each
$x=(\mbox{}_0x,...,\mbox{}_{n-1}x)\in \bigoplus _{j=0}^{n-1}
(Y\oplus Y{\bf i})$. Then a sequence $ \{ B(D(\mbox{}_k{\cal M})x):
~ k\in {\bf N} \} $ $~c$-converges to $B(D({\cal K})x)$, since an
operator $B$ is $c$-continuous. By Defintions 3 and 37 this means
that the limit exists $$\lim_{k\to \infty } \| {\hat
B}(\mbox{}_kM)x- {\hat B}(K)x \| = \lim_{k\to \infty }
\max_{j=0}^{n-1} \| B(D(\mbox{}_k{\cal M})x)(j) - B(D({\cal K})x)(j)
\| =0$$ and hence $ \| {\sf F} \| \le 1$, where ${\sf F}$ denotes
the Fourier transform operator on $L_q^{n,per}(l_{\infty } ({\bf Z
},Y))$ with values in $C_s({\bf S}^1,(L_q(Y\oplus Y{\bf i}))^n)$.

\par {\bf 39. Corollary.} {\it If a sequence $\{ B_p: p\in {\bf N}
\} $ of $n$-periodic operators converges to an $n$-periodic operator
$B$ relative to the norm on $L_q^{n,per}(l_{\infty } ({\bf Z },Y))$,
then a sequence of their Fourier transforms ${\sf F}(B_p)$ converges
to ${\sf F}(B)$ in $C_s({\bf S}^1,(L_q(Y\oplus Y{\bf i}))^n)$.}

\par {\bf 40. Corollary.} {\it If $B$ is an $n$-periodic operator
and ${\sf F}(B)={\hat B}$ its Fourier transform, then $\| B \| \ge
\sup_{M\in {\bf S}^1} \| {\hat B}(M) \| $.}

\par {\bf 41. Notation.} A family of all $n$-periodic operators
$B\in L_q^{n,per}(l_{\infty } ({\bf Z },Y))$ such that its Fourier
transform ${\sf F}(B)=\hat B$ has an absolutely converging Fourier
series $${\hat B}(M) = \sum_{k= -\infty }^{\infty }
\bigoplus_{j=0}^{n-1} M^{(k-1)n+j} \mbox{ }_{(k-1)n+j}B,$$ i.e. $
\sum_{k= -\infty }^{\infty }\max_{j=0}^{n-1} \| \mbox{}_{(k-1)n+j}B
\| <\infty $, will be denoted by $L^{n,1}_q(l_{\infty }({\bf
Z},Y))$.

\par {\bf 42. Lemma.} {\it Let $B\in L^{n,per}_q(l_{\infty }({\bf Z},Y))$
be an $n$-periodic operator and its Fourier transform ${\sf
F}(B)=\hat B$ has the form:
$${\hat B}(M) = \sum_{l=-\infty }^{\infty } M^{l} \mbox{}_{l}B ,$$
where $Y$ is a Banach space over the Cayley-Dickson algebra ${\cal
A}_v$, $~2\le v$. Then
$$(1)\quad B = \sum_{k=-\infty }^{\infty }\mbox{}_k\bar{B} S(kn),$$
where $\mbox{ }_k\bar{B} = \bigoplus_{j=0}^{n-1} \mbox{
}_{(k-1)n+j}B\in \bigoplus_{j=0}^{n-1} L_q^c(l_{\infty }({\bf
Z},Y))$ is an operator of (left) multiplication on stationary
operator valued sequence $\mbox{}_k\bar{B}\in \bigoplus_{j=0}^{n-1}
L_q(Y),$ moreover, \par $(2)$ $ \| \mbox{}_k\bar{B} \| = \sup_{0\le
j\le n-1} \| \mbox{}_{(k-1)n+j} B \| $ for each $k$.}
\par {\bf Proof.} In view of Lemma 7 the Banach spaces $\bigoplus_{j=0}^{n-1} L_q^c(l_{\infty }({\bf
Z},Y))$ and $L_q^c(l_{\infty }({\bf Z},Y))$ are isometrically
isomorphic, that follows from using the block form of matrices
$(B_{(k-1)n+j,(m-1)n+l})$ of operators $B$, where $k, m\in \bf Z$
and $j, l=0,...,n-1$, $~ n\ge 1$. From Lemma 5 it follows that an
operator $\sum_{k=-\infty }^{\infty }\mbox{}_{(k-1)n+j}B S(kn)$ is
$c$-continuous for each $j$, consequently, $\sum_{k=-\infty
}^{\infty }\mbox{}_k\bar{B} S(kn)$ is also $c$-continuous. A matrix
of the operator $B$ coincides with that of $\sum_{k=-\infty
}^{\infty }\mbox{}_k\bar{B} S(kn)$ by Lemma 35 and Definition 37.
Therefore, Formula $(1)$ is satisfied in accordance with Lemma 7.
The natural isometric embedding $Y\hookrightarrow Y\oplus Y{\bf i}$
induces isometric embeddings $L_q(Y)\hookrightarrow L_q(Y\oplus
Y{\bf i})$ and $L_q(l_{\infty }({\bf Z},Y))\hookrightarrow
L_q(l_{\infty }({\bf Z},Y\oplus Y{\bf i}))$ of normed spaces over
the Cayley-Dickson algebra ${\cal A}_v$.
\par From the definition of the operator norms on $L_q(Y\oplus Y{\bf i})$
and $\bigoplus_{j=0}^{n-1} L_q(Y\oplus Y{\bf i})$ (see Formulas
38$(1,2)$) and $L_q^c(l_{\infty }({\bf Z},Y\oplus Y{\bf i}))$
Equality $(2)$ follows, where \par $(3)$ $\| A \| = \sup_{0\le j\le
n-1} \| \mbox{}_jA \| $ \\ for each
$A=(\mbox{}_0A,...,\mbox{}_{n-1}A)\in \bigoplus_{j=0}^{n-1}
L_q^c(l_{\infty }({\bf Z},Y\oplus Y{\bf i}))$.
\par {\bf 43. Corollary.} {\it If $B, D\in L^{n,per}_q(l_{\infty }({\bf Z},Y))$,
then $$(1)\quad (Bx)(l) =\sum_{s=-\infty }^{\infty }\mbox{}_{s}B
x(s+l)$$ $$= \sum_{s=-\infty }^{\infty } \mbox{}_{s-l}B x(s)=:
(b\star x)(l)$$ and $(BD)(x) = b\star (d\star x)$ for each $x\in
l_{\infty }({\bf Z},Y)$, where $b=\{ \mbox{}_sB: ~ s\in {\bf Z} \}
\in l_1({\bf Z},L_q(Y))$, where $Y$ is a Banach space over the
Cayley-Dickson algebra ${\cal A}_v$. Particularly, as $0\le v\le 2$
the convolution is associative $b\star (d\star x)=(b\star d)\star
x$. }
\par {\bf Proof.} The equalities follow
$$(Bx)(l) =\sum_{k=-\infty }^{\infty
}\sum_{j=0}^{n-1}\mbox{}_{(k-1)n+j}B x((k-1)n+j+l)$$ $$=
\sum_{k=-\infty }^{\infty }\sum_{j=0}^{n-1}\mbox{}_{(k-1)n+j-l}B
x((k-1)n+j)$$ from Lemma 42. Putting $s=(k-1)n+j$ one gets Formula
$(1)$.
\par {\bf 44. Remark.} If $n=1$, the Fourier transform of an
operator valued function $b: {\bf Z}\to L_q(Y)$ with $b\in l_1({\bf
Z},L_q(Y))$ coinsides with the Fourier series for a mapping
$\breve{B}$.

\par {\bf 45. Proposition.} {\it Let $A$ and $B$ be two operators in
$L^{n,per}_q(l_{\infty }({\bf Z},Y))$, where $Y$ is a Banach space
over the Cayley-Dickson algebra ${\cal A}_v$, $~2\le v$. Then
$$(1)\quad \widehat{AB}(M) = \hat{A}(M) \hat{B}
(M)$$ for each $M\in {\bf S}^1$.}
\par {\bf Proof.} With $M\in {\bf S}^1$
we infer that
$$(2)\quad \mbox{}_m(AB) = \sum_{p=0}^m\mbox{
}_pA\mbox{ }_{m-p}B$$ and this implies Formula $(1)$, since
$\widehat{AB}(M)x=A(D({\cal M})B(D({\cal M})x))$ for each $x\in Y$
and $M\in {\bf S}^1$, since ${\bf i}i_j = i_j{\bf i}$for each $j$.

\par {\bf 46. Proposition.} {\it Let an operator
$B\in L^{n,per}_q(l_{\infty }({\bf Z},Y))$ be $n$-periodic, where
$Y$ is a Banach space over the Cayley-Dickson algebra ${\cal A}_v$,
$~2\le v$. Then an operator ${\hat B}(M)$ is invertible and $({\hat
B}(M))^{-1} x = \widehat{A^{-1}}(M)x$ for each $M\in {\bf S}^1$ and
$x\in Y$.}
\par {\bf Proof.} If $N\in {\bf S}^1$, then an algebra
$alg_{\bf R}(N,i_s)$ is associative for each $s\ge 0$, since
$N=N_0+N_1{\bf i}$ with $N_0, N_1\in \bf R$ and ${\bf i}i_s =
i_s{\bf i}$ for each $s\ge 0$. If $M\in {\bf S}^1$ and $x\in Y$, one
can take the algebra $alg_{\bf R}(M)$ which is either the real or
complex field. Therefore, $B(D({\cal M})B^{-1}(D({\cal M})x)) =
AA^{-1}x=x$ with $Ax = B(D({\cal M})x)$ by Proposition 45.

\par {\bf 47. Corollary.} {\it Let $B, D \in L^{n,per}_q(l_{\infty }({\bf
Z},Y))$ be $n$-periodic operators and let $({\hat B}(M))^{-1}x={\hat
D}(M)x$ for each $M\in {\bf S}^1$ and $x\in Y$, where $Y$ is a
Banach space over the Cayley-Dickson algebra ${\cal A}_v$, $~2\le
v$. Then $D=B^{-1}$.}
\par {\bf Proof.} This follows from Proposition 46, since the $\bf
R$ linear span $span_{\bf R} \{ y=Mx: ~ M\in {\bf S}^1, ~ x \in X \}
$ of such set of vectors is isomorphic with $X\oplus X{\bf i}$.

\par {\bf 48. Lemma.} {\it Let an $n$-periodic operator
$B\in L_q^{n,per}(l_{\infty }({\bf Z },Y))$ be uniformly
$c$-continuous, where $Y$ is a Banach space over the Cayley-Dickson
algebra ${\cal A}_v$, $~2\le v$. Then its Fourier transform $\hat B$
is uniformly $c$-continous, ${\hat B}\in L_q^{uc}(l_{\infty }({\bf
Z},Y\oplus Y{\bf i}))$.}
\par {\bf Proof.} From the conditions of this lemma it follows, that
the mapping $\breve{B}: {\bf S}^1\ni M\mapsto D({\cal M}) B D({\cal
M}^*)$ is continuous from ${\bf S}^1$ into $L_q^{n,per}(l_{\infty
}({\bf Z },Y\oplus Y{\bf i}))$. Up to an $\bf R$-linear continuous
algebraic automorphism of the Cayley-Dickson algebra ${\cal A}_v$
and the corresponding automorphism of a Banach space $Y$, the
Fourier series
$$(1)\quad \breve{B}(M) \sim \sum_{k= -\infty }^{\infty }
M^k\mbox{ }_k\breve{B}$$ exists by Lemma 31. This series converges
to $\breve{B}(M)$ by Cezaro, that is $$(2)\quad \breve{B}(M)
=\lim_{m\to \infty } \sum_{k= -m}^m(1- \frac{|k|}{m+1}) M^k\mbox{
}_k\breve{B},$$ since $(1- \frac{|k|}{m+1})\in \bf R$ while the real
field is the center of the Cayley-Dickson algebra ${\cal A}_v$.
Particularly, for $M=1$ one has $M^k=1$ and $\breve{B}(1)=B$ (see
also \S 20.2(743) \cite{fihteng}).
\par The operator $B$ is $n$-periodic, so consider its Fourier
transform and get
$$(3)\quad {\hat B}(M) =\lim_{m\to \infty } \sum_{k= -m}^m
(1- \frac{|k|}{m+1}) \mbox{ }_k\hat{B}(M).$$ In view of Lemma 31
each operator $\mbox{ }_k\breve{B}S(-k)$ is diagonal. Since $B$ is
$n$-periodic, this implies that every operator
$\bigoplus_{j=0}^{n-1}\mbox{ }_{(k-1)n+j}\breve{B}S(-(k-1)n-j)$ is
$n$-periodic as well. Therefore, $\mbox{}_k\hat{B}(M)x = (\mbox{
}_k\breve{B}(D({\cal M})x))(0)=(M^kI)\mbox{ }_kBx$ for each vector
$x\in Y\oplus Y{\bf i}$, since
$M^k(M^{-k}x_li_l)=M^k(M^{-k}i_lx_l)=x_li_l$ for each $l\ge 0$ and
$x_l\in Y_l$, consequently, $\mbox{}_k\hat{B}(M)= (M^kI)\mbox{ }_kB$
and hence
$$(4)\quad {\hat B}(M) =\lim_{m\to \infty } \sum_{k= -m}^m
(1- \frac{|k|}{m+1}) (M^kI)\mbox{ }_kB,$$ where
$(\mbox{}_k\breve{B})_{s,p} = \mbox{}_kB$ for each $s-p=k$, $~ s, p
\in \bf Z$.

\par {\bf 49. Notation.} Let $\sf P$ be a Banach algebra over the
Cayley-Dickson algebra ${\cal A}_v$ with $2\le v$. We denote by
$F({\bf S}^1,{\sf P})$ a Banach space of all continuous functions
$f: {\bf S}^1\to \sf P$ with absolutely converging Fourier series
$$(1)\quad f(M) = \sum_{k= - \infty }^{\infty } M^k \mbox{ }_kf$$
relative to the norm: $$(2)\quad \| f \| := \sum_{k= - \infty
}^{\infty } \| \mbox{ }_kf \|,$$ where $\mbox{}_kf\in \sf P$ for
each $k\in \bf Z$.

\par {\bf 50. Corollary.} {\it Let $B\in L_q^{n,per}(l_{\infty }
({\bf Z},Y))$, where $Y$ is a Banach space over the Cayley-Dickson
algebra ${\cal A}_v$, $~2\le v$. Then the following conditions are
equivalent:
\par $(1)$ $\hat{B} \in F({\bf S}^1,L_q(Y\oplus Y{\bf i}))$ and
\par $(2)$ $\breve{B}\in F({\bf S}^1,L_q^c(l_{\infty }
({\bf Z},Y\oplus Y{\bf i}))$.}
\par {\bf Proof.} This follows from Lemma 48, since $\|
\mbox{}_k\breve{B} \| = \| \mbox{}_kB \| $ for each $k\in \bf Z$.
Indeed, generally $\| \mbox{}_k\breve{B} x\| = \| \mbox{}_kB x\| $,
since $ |ab|\le |a| |b|$ for each Cayley-Dickson numbers $a, b\in
{\cal A}_v$ and $ \| a x \| \le |a| \| x \| $ for each $a\in {\cal
A}_v$ and $x\in Y$ (see \S I.2.1 \cite{ludopalglamb}). In
particular, if $x\in Y_0$ or $x\in Y_l$, then $ \| a x \| =|a| \| x
\| $.

\par {\bf 51. Theorem.} {\it The algebra $L_q^{n,per}(l_{\infty }
({\bf Z},Y))$ is the subalgebra of the algebra $L_q^{uc}(l_{\infty }
({\bf Z},Y))$, where $Y$ is a Banach space over the Cayley-Dickson
algebra ${\cal A}_v$, $~2\le v$.}
\par {\bf Proof.} If $B\in L_q^{n,per}(l_{\infty }
({\bf Z},Y))$, then by Lemma 42 we have $$B = \sum_{k=-\infty
}^{\infty }\mbox{}_k\bar{B} S(kn),$$ where $\mbox{ }_k\bar{B} =
\bigoplus_{j=0}^{n-1} \mbox{ }_{(k-1)n+j}B \in \bigoplus_{j=0}^{n-1}
L_q^c(l_{\infty }({\bf Z},Y))$ is an operator of (left)
multiplication on stationary operator valued sequence
$\mbox{}_k\tilde{B} \in \bigoplus_{j=0}^{n-1} L_q(Y),$ moreover,
\par $ \| \mbox{}_k\bar{B} \| = \| \mbox{}_k\tilde{B} \| $ for each $k$. In
view of Proposition 24 an $1$-ribbon operator $\mbox{}_k\bar{B}$ is
uniformly $c$-continuous. On the other hand, each shift operator
$S(n)$ is uniformly $c$-continuous. The algebra $L_q^{uc}(l_{\infty
} ({\bf Z},Y))$ is complete as the uniform space. Therefore, the
operator $B$ is uniformly $c$-continuous.

\section{Fourier transform on algebras and spectra}

\par {\bf 52. Definitions.} Let $\sf G$ be a quasi-group,
i.e. a set with one binary operation (multiplication) so that
\par $(1)$ there exists a unit element $e$ so that $eb=be=b$;
\par $(2)$ each element $b$ has an inverse $b^{-1}$, i.e.
$b^{-1}b=bb^{-1}=e$;
\par $(3)$ a multiplication is alternative $(aa)b=a(ab)$ and
$b(aa)=(ba)a$ and \par $(4)$ $a^{-1}(ab)=b$ and $(ba)a^{-1}=b$ for
each $a, b \in {\sf G}$.
\par Let $R_a$ be a Banach algebra over the real field $\bf R$ for each $a\in {\sf G}$ such that
$R_a$ is isomorphic with $R_b$ for all $a, b\in {\sf G}$.  Put
$$(5)\quad R= \{ B: ~ B\in  \bigoplus_{a\in {\sf G}} a R_a, ~ \| B \| <\infty \}
,$$  $$ ~ {\sf G}_{\bf R} =: \{ x: x\in \bigoplus_{a\in {\sf G}}
a{\bf R}, ~ |x|<\infty \} $$ is a quasi-group ring over the real
field so that $a\beta =\beta a$ for each $a\in {\sf G}$ and $\beta
\in {\bf R}$,
$$(6)\quad |x|^2=\sum_{a\in {\sf G}} |x_a|^2\mbox{ for }
x=\sum_{a\in {\sf G}} x_aa$$ with $x_a\in {\bf R}$ for each $a\in
{\sf G}$;
$$(7)\quad \| A \| ^2 := \sum_{b\in {\sf G}} \| A_b \| ^2$$ for $A=\sum_{b\in {\sf G}} A_bb$
with $A_b\in R_b$ for each $b\in {\sf G}$, \par $(8)$ $bA_a=A_ab$
for each $A_a\in R_a$ and $a, b \in {\sf G}$. Suppose that $R_a$
\par $(9)$ contains a unit element $I$ and that \par $(10)$ $ \| I \| =1$ and
\par $(11)$ $ \| AB \| \le \| A \| \| B \| $ for each $A, B\in R_a$.

\par Denote by $L(R)$ the Banach space of all strongly integrable
functions $f: {\bf R}\to R$ supplied with the norm
$$(12)\quad \| f \| := \int_{-\infty }^{\infty } \| f(t) \| dt <\infty
.$$ Henceforward, we suppose that an algebra $R$
\par $(13)$ is alternative $(AA)B=A(AB)$ and
$B(AA)=(BA)A$ for each $A, B \in R$ and \par $(14)$  if $A$ is left
invertible, then also $A^{-1}(AB)=B$, if $A$ is right invertible
$(BA)A^{-1}=B$ for every $A, B \in R$.
\par The alternativity implies that ${\sf G}$ and $R$ are
power-associative that is $b^mb^n=b^{n+m}$ and $B^mB^n=B^{n+m}$ for
each $b\in {\sf G}$ and $B\in R$ and natural numbers $n, m$, where
$b^n=b(b(...(bb)...))$ denotes the $n$-fold product, $b^0=e$ for
$b\ne 0$, $~B^0=I$ for $B\ne 0$.
\par We consider their complexifications ${\sf G}_{{\bf C}_{\bf i}}
: = {\sf G}_{\bf R}\oplus {\sf G}_{\bf R}{\bf i}$ and $R_{{\bf
C}_{\bf i}} := R\oplus R{\bf i}$, where
 ${\bf C}_{\bf i} = {\bf R}\oplus {\bf R}{\bf i}$ with $a{\bf
i}={\bf i}a\in {\sf G}_{{\bf C}_{\bf i}}$, \par $(15)$ $~ | a+b{\bf
i}|^2 = |a|^2+|b|^2$ for each $a, b\in {\sf G}_{\bf R}$ and \par
$(16)$ $\| A+B{\bf i} \| ^2 = \| A \| ^2 + \| B \| ^2$ for every $A,
B\in R$. \par Analogously a Banach space $L(R_{{\bf C}_{\bf i}})$ is
defined with $f: {\bf R}\to R_{{\bf C}_{\bf i}}$.

\par {\bf 53. Lemma.} {\it Let $R$ be an algebra as in \S 52 and let an element
$A\in R$ be of norm $ \| A \| <1$, then the series
$C=I-A+A^2-A^3+...$ is absolutely convergent and \par $(1)$
$C(I+A)=(I+A)C=I$.}
\par {\bf Proof.} From Formulas 52$(7,11,16)$ it follows that
$\| A^n \| \le \| A \| ^n$ for each natural number $n$,
consequently, the sequence of partial sums $S_n :=
I-A+A^2-A^3+...+(-1)^nA^n$ converges in $R$. A Banach algebra $R$ is
power-associative and this implies Formula $(1)$.
\par {\bf 54. Lemma.} {\it Let $R$ be an algebra as in \S 52 and let
an element $A\in R$ have a left inverse $Q$. If $B\in R$ is an
element such that $ \| B \| \| Q \| <1$, then $(A+B)$ has a left
inverse $C$ so that
\par $(2)$ $C=Q(I-BQ+(BQ)^2-(BQ)^3+...)$.}
\par {\bf Proof.} A Banach algebra $R$ satisfies conditions 52$(13,14)$, hence
$(A+B)=(I+BQ)A$, since $QA=I$. The alternativity $(13)$ implies the
Moufang identities in the algebra $R$: \par $(M1)$  $(XYX)Z =
X(Y(XZ))$, \par $(M2)$ $Z(XYX) = ((ZX)Y)X$,
\par $(M3)$ $(XY)(ZX) = X(YZ)X$ for each $X, Y, Z\in R$.
\par From Lemma 53 and Formulas $(M1,M2)$ it follows that
$C(A+B)=(Q(I-BQ+(BQ)^2-...))((I+BQ)A)=(Q-Q(BQ)+Q(BQ)^2-...)(A+B)=I$.
\par {\bf 55. Corollary.} {\it The set $U_l$ of all left invertible
elements $A\in R$ is an open subset in $R$.}
\par {\bf 56. Notation.} Denote by $R'$ the algebra over ${\bf R}$ of all periodic functions
$x: [0,2\pi ]\to R_{{\bf C}_{\bf i}}$ of the form
$$(1)\quad x(t) = \sum_{n= - \infty }^{\infty } a_n e^{nt\bf i}$$
with coefficients $a_n\in R$ such that
\par $$(2)\quad \sum_n \| a_n \| <\infty $$ with point-wise addition and
multiplication of functions
\par $$(3)\quad (x+y)(t)=x(t)+y(t), \quad (xy)(t)=x(t)y(t)\mbox{  for each  } t\in [0,2\pi ].$$
\par {\bf 57. Lemma.} {\it Suppose that $x\in R'$ and $x(0)$ has a left inverse
in $R$, then there exists an element $$(1)\quad y(t)=\sum_{n= -
\infty }^{\infty } c_ne^{nt\bf i}\in R'$$ such that
\par $(2)$ $c_0$ has a left inverse $q_0$ in $R$ and \par $(3)$ $\|
q_0 \| \sum_{n=1}^{\infty } \| c_n +c_{-n} \| <1$ and \par $(4)$
there exists $\epsilon >0$ so that $y(t) = x(t)$ for each $t\in
(-\epsilon ,\epsilon )$.}
\par {\bf Proof.} Consider the following function given piecewise
$w_{\epsilon }(t)=1$ for $|t|<\epsilon $, $~w_{\epsilon } (t) = 2 -
|t|/\epsilon $ for $\epsilon \le |t|<2\epsilon $, $~w_{\epsilon
}(t)=0$ for $2\epsilon \le |t|$, where $0<\epsilon \le \pi /2 $.
Then one defines the function $$ y_{\epsilon }(t) = w_{\epsilon
}(t)x(t) + [1-w_{\epsilon }(t)]x(0)= \sum_{n= -\infty }^{\infty }
b_n(\epsilon )e^{nt\bf i} .$$ This function satisfies Condition
$(4).$ It has the Fourier series with coefficients $b_n=b_n(\epsilon
)$:
$$b_n=\frac{3\epsilon }{2\pi } a_n +\sum_{k=1}^{\infty }
\frac{a_{n-k}+ a_{n+k}}{\pi k^2\epsilon } (\cos (\epsilon k) - \cos
(2\epsilon k)) - \sum_{k=-\infty }^{\infty } a_k \frac{\cos
(\epsilon n)-\cos (2\epsilon n)}{\pi n^2\epsilon }$$ for $n\ne 0$
and
$$b_0= a_0 + \sum_{k=1}^{\infty } (a_{-k}+a_k)[1+\frac{\cos (\epsilon
k) - \cos (2\epsilon k)}{\pi k^2\epsilon } - \frac{3\epsilon }{2\pi
}].$$  Therefore, $$\lim_{\epsilon \downarrow 0} \| b_0 \| = \|
\sum_{k= - \infty }^{\infty } a_k \| = \| y(0) \| >0\mbox{  and}$$
$$\sum_{k=1}^{\infty } [\| b_k \| + \| b_{-k} \| ]\le \sum_{k=-\infty }^{\infty
}\| a_k \| A_k,$$ where a positive number $\delta>0$ exists such
that $0\le A_k <\epsilon ^{1/2} [2|k|C+9/\pi ]$ for each $0<\epsilon
<\delta $ and every $k\in \bf Z$, where $C=const >0$. Thus a
positive number $\epsilon _0>0$ exists so that $$\| b_0\| >
\sum_{k=1}^{\infty }[\| b_k \| + \| b_{-k} \| ]$$ for each
$0<\epsilon <\epsilon _0$ (see also \cite{annmbochph42,wienannm32}).
From Lemma 54 statements $(2,3)$ of this lemma follow.
\par {\bf 58. Corollary.} {\it If $y\in R'$ and $y(t)$ satisfies
Properties $(2,3)$ of Lemma 57, then $y$ has a left inverse $z$ in
$R'$.}
\par {\bf 59. Theorem.} {\it If $R'$ is an algebra of periodic
functions as in \S 56. Then $x(t)$ has a left inverse in $R'$ if
$x(t_0)$ has a left inverse in $R$ for each $t_0$.}
\par {\bf Proof.} In view of Lemma 57 and Corollary 58 for each $\tau \in [0,2\pi ]$
a positive number $\epsilon >0$ and an element $y_{\tau }\in R'$
exist such that $y_{\tau }(t)x(t)=I$ for each $t\in (\tau - \epsilon
,\tau +\epsilon )$. The segment $[ - \pi , \pi ]$ is compact, that
is, each its open covering has a finite subcovering, consequently, a
finite number of functions $y_{\tau }$ induces a function $y\in R'$
so that $y(t)x(t)=I$ for each $t$.
\par {\bf 60. Lemma.} {\it Suppose that $ - \pi <\alpha <a<b< \beta <\pi $
and $x_1, x_2\in R'$ and $x_2(t)$ has a left inverse for each $t\in
(\alpha ,\beta )$ and $x_1(t)$ vanishes for every $- \pi \le t<a$
and $b<t\le \pi $. Then an element $x_3\in R'$ exists vanishing on
$[-\pi ,\alpha )\cup (\beta ,\pi ]$ such that
\par $(1)$ $x_1(t)=x_3(t)x_2(t)$ for each $t\in [-\pi ,\pi ]$. }
\par {\bf Proof.} From Lemma 57 and Corollary 58 it follows, that
to any $\tau \in [a,b]$ a positive number $\epsilon >0$ and an
element $y_{\tau } \in R'$ correspond such that $y_{\tau
}(t)x_2(t)=I$ for each $t\in (\tau - \epsilon, \tau +\epsilon )$. As
in \S 59 one gets that an element $z\in R'$ exists such that
$z(t)x_2(t)=I$ for every $t\in [a,b]$. Put $x_3(t)=x_1(t)z(t)$,
consequently, $x_3(t)=0$ for each $t\in [-\pi ,\alpha ]\cup [\beta
,\pi ]$ and $x_3\in R'$ and hence Assertion $(1)$ is valid, since
the algebra $R$ satisfies Conditions 52$(13,14)$ and $e^{nt{\bf i}}$
commutes with $R$ and $(a_ne^{nt{\bf i}})(b_ke^{kt{\bf i}})=
(a_nb_k) e^{(n+k)t{\bf i}}$ for each $a_n, b_k\in R$ and $n, k \in
{\bf Z}$.

\par {\bf 61. Lemma.} {\it Suppose that $x(t)$ is strongly integrable on $(-\pi ,\pi )$ function
with values in $R$ and vanishes on $(-\pi , -\pi +\epsilon )\cup
(\pi - \epsilon ,\pi )$ with $0<\epsilon <\pi /2$ and $$(1)\quad
f(t) = \frac{1}{2\pi } \int_{-\pi }^{\pi } x(\tau )e^{-\tau t{\bf
i}}d\tau ,$$ $$(2)\quad a_n=\frac{1}{2\pi } \int_{-\infty }^{\infty
} x(\tau ) e^{-n\tau {\bf i}} d\tau ,$$ then $$(3)\quad
\int_{-\infty }^{\infty } \| f(t) \| dt <\infty $$ if and only if
\par $$(4)\quad \sum_{n=-\infty }^{\infty } \| a_n \| <\infty .$$}
\par {\bf Proof.} Consider a positive number $0<\delta <\pi /2$ so that
$x(t)$ vanishes on $[-\pi ,-\pi +2\delta )\cup (\pi -2\delta ,\pi
)$. Put $\phi (t)=1$ for $|t|<\pi - \delta $, $\phi (t)=\frac{\pi -
|t|}{\delta }$ for $|\pi - \delta | \le |t|<\pi $, $\phi (t)=0$ for
$\pi \le |t|$. If
$$(5) \quad x(t) = \sum_{n= -\infty }^{\infty } a_n e^{nt{\bf i}},$$
then
$$(6)\quad x(t)=x(t)\phi (t) = \frac{1}{\pi } \int_{-\infty
}^{\infty } e^{\tau t{\bf i}} [\sum_{n=-\infty }^{\infty } a_n
\frac{\cos (\tau +n) (\pi -\delta ) - \cos (\tau +n) \pi } {(\tau
+n)^2\epsilon } ]d\tau ,$$ since this integral and this sum are
absolutely convergent. Therefore, the function
$$(7)\quad f(t) := \sqrt{\frac{2}{\pi } } [\sum_{n=-\infty }^{\infty } a_n \frac{\cos (\tau +n) (\pi
-\delta ) - \cos (\tau +n)\pi }{(\tau +n)^2\epsilon }]$$ satisfies
Conditions $(1,3)$, if $(4)$ is fulfilled.
\par Vise versa, Condition $(3)$ implies that
$$(8)\quad \sum_{n= - \infty }^{\infty } \| \int_{n-1/2}^{n+1/2} f(t) dt \|  <\infty ,$$
consequently,
$$(9) \quad \int_{n-1/2}^{n+1/2} f(t) dt = \frac{1}{\sqrt{2\pi }}
\int_{-\pi }^{\pi } x(t)\frac{2\sin (t/2)}{t} e^{nt{\bf i}} dt .$$
Thus the Fourier series of the function $x(t)\frac{2\sin (t/2)}{t}$
converges absolutely. Moreover, the Fourier series of the mapping
$\frac{t}{2\sin (t/2)} \phi (t)$ also is absolutely convergent. Thus
the Fourier series of $x(t)=[x(t)\frac{2\sin
(t/2)}{t}][\frac{t}{2\sin (t/2)} \phi (t)]$ is absolutely
convergent, since $R'$ is an algebra over the real field $\bf R$ and
$\bf i$ commutes with each $y\in R'$.
\par {\bf 62. Corollary.} {\it Let $g$ and $f\in L(R)$, let also
$$(1)\quad x_1(t) = \int_{-\infty }^{\infty } g(\tau )e^{-\tau t{\bf i}} d\tau $$
vanish outside some interval $(a,b)\subset (-\pi ,\pi )$. Suppose
that
$$(2)\quad x_2(t) = \int_{-\infty }^{\infty } f(\tau ) e^{-\tau t{\bf i}} d\tau $$
is zero outside an interval $(\alpha ,\beta )\subset (-\pi ,\pi )$
with $\alpha <a$ and $b<\beta $ and $x_2(t)$ has a left inverse for
each $\alpha <t<\beta $. Then an element $y\in L(R)$ exists so that
$$(3)\quad g(t)= \int_{-\infty }^{\infty } y(\tau ) f(t-\tau )d\tau
.$$ }
\par This follows immediately from Lemmas 60 and 61.

\par {\bf 63. Lemma.} {\it If $f\in L(R)$, then $$(1)\quad \lim_{\epsilon \downarrow 0}
\int_{-\infty }^{\infty } \| f(t+\epsilon ) - f(t) \| dt =0.$$}
\par {\bf Proof.} The Lebesgue measure is
$\sigma $-finite and $\sigma $-additive, so the statement of this
theorem for step functions is evident. In $L(R)$ the set of step
functions $$g(t)=\sum_{k=1}^n a_k \chi _{B_k}(t)$$ is dense with
$a_k\in R$, $~B_k\in {\cal B}({\bf R})$, $~n\in \bf N$, where ${\cal
B}({\bf R})$ denotes the $\sigma $-algebra of all Borel subsets of
$\bf R$.
\par {\bf 64. Lemma.} {\it Suppose that $f\in L(R)$ and $h\in L(R)$ so that $supp (h)\subset
(-\epsilon ,\epsilon )$ for some positive number $0<\epsilon <\infty
$. Then $$(1)\quad \int_{-\infty }^{\infty } \| f(t) \int_{-\infty
}^{\infty } h(\tau )d\tau  - \int_{- \infty }^{\infty } f(t+\tau
)h(\tau ) d\tau \| dt $$ $$\le [\int_{-\infty }^{\infty } \| h(\tau
) \| d\tau ] \sup_{|u|\le \epsilon } \int_{-\infty }^{\infty } \|
f(t+u)-f(t) \| dt .$$ }
\par {\bf Proof.} This follows from Fubini's theorem
$$\int_{-\infty }^{\infty } \| f(t) \int_{-\infty
}^{\infty } h(\tau )d\tau  - \int_{- \infty }^{\infty } f(t+\tau
)h(\tau ) d\tau \| dt $$ $$\le \int_{-\infty }^{\infty } \| f(t)
 -  f(t+\tau ) \| \| h(\tau ) \| d\tau dt $$ $$\le \int_{-\infty }^{\infty } \| h(\tau ) \| d\tau
\sup_{|u|\le \epsilon } \int_{-\infty }^{\infty } \| f(t+u)-f(t) \|
dt .$$
\par {\bf 65. Lemma.} {\it If $f\in L(R)$, then
$$(1)\quad \lim_{n\to \infty } \int_{-\infty }^{\infty } \| f(t) -
\frac{1}{\pi n} \int_{-\infty }^{\infty } f(t+\tau )\frac{\sin
^2(n\tau )}{\tau ^2} d\tau \| dt =0.$$ }
\par {\bf Proof.} Put $h_0(t)=\frac{\sin
^2(nt )}{t ^2}=h_1(t)+h_2(t)$ with $h_1(t)=h_0(t)[1-|t|\sqrt{n}]$
for $|t|\le n^{-1/2}$, while $h_1(t)=0$ for $|t|>n^{-1/2}$. An
application of Lemmas 63 and 64 leads to
$$(2)\quad \lim_{n\to \infty } \int_{-\infty }^{\infty } \| f(t) -
\frac{1}{\pi n} \int_{-\infty }^{\infty } f(t+\tau )h_1(\tau ) d\tau
\| dt =0$$  and
$$(3)\quad \lim_{n\to \infty } \int_{-\infty }^{\infty } \frac{1}{\pi n} \|
f(t+\tau )h_2(\tau ) d\tau \| dt $$
$$\le [\int_{-\infty }^{\infty } \| f(t) \| dt] \lim_{n\to \infty }
\frac{1}{\pi n} \int_{- \infty }^{\infty }h_2(\tau )d\tau .$$ On the
other hand,
$$(4)\quad \frac{1}{\pi n} \int_{-\infty }^{\infty } h_2(t)dt=$$
$$\frac{1}{\pi n} [\int_{-\infty }^{\infty } \frac{\sin ^2(nt )}{t
^2}dt - \int_{-n^{-1/2}}^{n^{-1/2}} (1-|t|\sqrt{n}] \frac{\sin ^2(nt
)}{t ^2} dt ]$$  $$=\frac{2}{\pi } \int_{n^{1/2}}^{\infty }
\frac{\sin ^2(t )}{t ^2} dt + \frac{2}{\pi \sqrt{n}}
\int_0^{\sqrt{n}} \frac{\sin ^2(t )}{t ^2} dt = O (n^{-1/2} \ln n)
,$$ consequently,
$$(5)\quad \lim_{n\to \infty } \int_{-\infty }^{\infty } \|
\frac{1}{\pi n} \int_{-\infty }^{\infty } f(t+\tau )h_2(\tau )d\tau
\| dt =0.$$
\par {\bf 66. Theorem.} {\it Let $f\in L(R)$ and let the Fourier transform
$$(1)\quad x(\tau ) = \int_{-\infty }^{\infty } f(t)e^{\tau t{\bf i}} dt $$
have a left inverse in $R$ for each $\tau \in [-\pi ,\pi ]$. Then
the $\bf R$-linear combinations $$(2)\quad \sum_n b_n f(t-\tau
_n)\mbox{
 with  } b_n\in R$$ are dense in $L(R)$, where $\tau _n\in [-\pi ,\pi )$.}
\par {\bf Proof.} Lemma 65 means that a function $$(3)\quad f_{\delta }(t) =
\frac{1}{\pi n} \int_{-\infty }^{\infty  } f(t+\tau ) \frac{\sin
^2(n\tau )}{\tau ^2}  d\tau $$ exists so that
$$(4)\quad \int_{-\infty }^{\infty } \| f(t) - f_{\delta }(t) \| dt
<\delta ,$$ where $0<\delta $. Consider its Fourier transform:
$$(5)\quad h_1(u) := \frac{1}{\sqrt{2\pi }} \int_{ - \infty }^{\infty }
e^{ut{\bf i}} [\frac{1}{\pi n} \int_{-\infty }^{\infty } f(t+\tau )
\frac{\sin ^2(n\tau )}{\tau ^2}d\tau ]dt$$
$$\frac{1}{\sqrt{2\pi }} \int_{- \infty }^{\infty } f(t)e^{tu{\bf i}} [\frac{1}{\pi n}
\int_{-\infty }^{\infty }\frac{\sin ^2(n\tau )}{\tau ^2}e^{-u\tau
{\bf i}} d\tau ]dt,\mbox{  consequently,}$$ $$h_1(u)
=(1-\frac{|u|}{2n})\frac{1}{\sqrt{2\pi }} \int_{- \infty }^{\infty }
f(t)e^{ut{\bf i}} dt\mbox{  when  } |u|<2n$$ and $h_1(u)=0$ for
$|u|\ge 2n$. Analogously the Fourier transform
$$h_2(u) := \frac{1}{\sqrt{2\pi }} \int_{ - \infty }^{\infty }
e^{ut{\bf i}} [\frac{1}{2\pi n} \int_{-\infty }^{\infty } f(t+\tau )
\frac{\sin ^2(2n\tau )}{\tau ^2}d\tau ]dt$$ vanishes for each $u$ so
that $|u|>4n$. The Fourier series of $h_1$ and $h_2$ over $(-8n,8n)$
converge absolutely by Lemma 61. Then one can write
$h_1(u)=h_2(u)h_3(u)$, where $$h_3(u)=\int_{- \infty }^{\infty }
\psi (t)e^{tu{\bf i}} dt $$ with $\psi \in L(R)$, since the algebra
$R$ is alternative and $e^{tu{\bf i}}$ commutes with any $y\in R'$
for each real numbers $t$ and $u$. Therefore, we deduce that
$$(6)\quad \int_{-\infty }^{\infty } f_{\delta }(t)e^{tu{\bf i}}
dt = $$
$$\int_{-\infty }^{\infty } \frac{e^{tu{\bf i}}}{2\pi n}
\{ [\int_{-\infty }^{\infty } f(t+\tau ) \frac{\sin ^2(2n\tau
)}{\tau ^2} ] \int_{-\infty }^{\infty } \psi (x)e^{ux{\bf i}} dx \}
d\tau ]dt \mbox{,  consequently,}$$
$$\int_{-\infty }^{\infty } e^{tu{\bf i}}  [ f_{\delta }(t) - \frac{1}{2\pi n}
\int_{-\infty }^{\infty } f(t+ x) \{ \int_{-\infty }^{\infty }
\frac{\sin ^2(2n\tau )}{\tau ^2} ]  \psi (\tau - x) d\tau \} dx ] dt
=0$$ for each $u$, hence
$$(7)\quad f_{\delta }(t) = \frac{1}{2\pi n} \int_{-\infty }^{\infty } f(t+ x) \Phi (x) dx\mbox{, where}$$
$$\Phi (x) = \int_{-\infty }^{\infty } \frac{\sin ^2(2n\tau )}{\tau ^2} ]  \psi
(\tau - x) d\tau $$ is absolutely integrable. Lemmas 63 and 64 imply
that
$$(8)\quad \lim_{n\to \infty } \int_{-\infty }^{\infty } \| \int_{-\infty }^{\infty }
f(t+x) \Phi (x)dx- \sum_{k=-n^2}^{n^2-1}f(t+ \frac{k}{n})
\int_{k/n}^{(k+1)/n} \Phi (x)dx \| dt=0.$$ From Formulas $(5,6)$ and
Lemma 65 the assertion of this theorem follows.
\par {\bf 67. Lemma.} {\it Let $R$ be an algebra with unit (see \S 52) and let $\cal I$
be a maximal left ideal. Suppose that $X$ is an additive group of
all cosets $R/{\cal I}$ and $H$ is an algebra of homomorphisms of
$X$ onto itself produced by multiplying by elements of $R$ from the
left. Then $X$ is irreducible relative to $H$.} \par {\bf Proof.}
Consider the quotient mapping $\theta : R\to R/{\cal I}$ (see also
\S I.2.39 \cite{ludopalglamb}). A space $X$ is $\bf R$-linear, since
$R$ is an algebra over $\bf R$. Therefore, $\theta (R)a =:V_a $ is
an $\bf R$-linear space for each nonzero element $a\in R\setminus \{
0 \} $. Put $S_a=\theta ^{-1}(V_a)$. Evidently, $S_a$ is a left
ideal in $R$ and ${\cal I}\subset S_a$, $~S_a\ne {\cal I}$. The
ideal $\cal I$ is maximal, hence $S_a=R$, consequently, $V_a=X$ and
$Ha=X$ for each nonzero element $a\in R\setminus \{ 0 \} $.
\par {\bf 68. Lemma.} {\it Let $R$, ${\cal I}$, $X$ and $H$ be the
same as in Lemma 67. Suppose that for a marked element $x\in R$ and
each maximal left ideal ${\cal I}$ the corresponding element $\theta
(x)$ is left invertible in $H$. Then $x$ is left invertible in $R$.}
\par {\bf Proof.} An algebra $R$ has the unit element $I\in R$. Each
left ideal is contained in a maximal left ideal. Therefore, an
element $x\in R$ is left invertible if and only if this element $x$
is not contained in a maximal left ideal. On  the other hand, if
$yx=I$, then $\theta (y)\theta (x)=\theta (I).$ Take cosets $b_I$
and $b_0$ in $R/\cal I$ so that $I\in b_I$ and $0\in b_0$. Then
$\theta (y)\theta (x) b_I = \theta (I) b_I$. But $\theta
(I)b_I=b_I$, since $II=I$. If $x\in b_0=I$, then $\theta (x)b_I=b_0$
and $\theta (y)\theta (x)b_I=b_0$, since $x\in \theta (x)b_I$.
\par {\bf 69. Lemma.} {\it A maximal left or right ideal $\cal I$ in
$R$ is closed.}
\par {\bf Proof.} If $\cal I$ is not closed, then its closure $cl_R({\cal
I})$ in $R$ contains an ideal $\cal I$. On the other hand,
$cl_R({\cal I})$ is a left or right ideal in $R$ respectively, since
$R$ is a topological algebra. Therefore, $cl_R({\cal I})={\cal I}$,
since $\cal I$ is maximal.

\par {\bf 70. Lemma.} {\it If $R$ is a Banach algebra, and ${\cal I}$, $X$ and $H$ have the
same meaning as in Lemma 67, then $X$ is a Banach space, $H$ is a
normed algebra with norm $ \| \theta (z) \| $ on $H$ so that $\|
\theta (x) \|_H \le \| x \|_R $ for each $x\in R$. If moreover $R$
is a Hilbert algebra over either the quaternion skew field or the
octonion algebra ${\cal A}_v$ with $2\le v \le 3$, then $X$ is a
Hilbert space over ${\cal A}_v$.}
\par {\bf Proof.} The quotient algebra $R/{\cal I}$ is supplied with
the quotient norm $ \| \theta (x) \|_X = \inf _{z\in \theta (x)} \|
z \|_R $ (see also \S I.2.39 \cite{ludopalglamb}). Therefore, $$ \|
\theta (x) \|_H = \sup_{b\in X} \| \theta (x) b \|_X / \| b \|_X \le
\| x \|_R .$$
\par If $R$ is a Hilbert algebra over ${\cal A}_v$, then $ \| x
\|_R =\sqrt{<x;x>}$, where a scalar ${\cal A}_v$-valued product
$<x;y>$ on $R$ satisfies conditions of \S I.2.3 \cite{ludopalglamb}.
From the parallelogram identity and the polarization formula one
gets that $ \| \theta (x) \|_X$ induces an ${\cal A}_v$-valued
scalar product on $X$ (see Formulas I.2.3$(1-3,SP)$
\cite{ludopalglamb}).

\par {\bf 71. Notation.} Let $R$ be a quasi-commutative $C^*$-algebra
over either the quaternion skew field or the octonion algebra ${\cal
A}_r$, $~2\le r \le 3$, satisfying conditions of \S 52. Let also $F$
be a normed algebra of functions either $f, g : \Lambda \to {\cal
A}_r$ or $({\cal A}_r)_{{\bf C}_{\bf i}}$ with point-wise
multiplication $f(t)g(t)$ and addition $f(t)+g(t)$ of functions and
$\phi : F\to {\cal A}_r$ or $({\cal A}_r)_{{\bf C}_{\bf i}}$
respectively be a continuous $\bf R$ homogeneous additive
multiplicative homomorphism, $\phi (fg)=\phi (f)\phi (g)$. Suppose
that the unit function $h(t)=1$ for each $t\in \Lambda $ belongs to
$F$, also $~R'$ is a family of functions $x: \Lambda \to R$ or $x:
\Lambda \to R_{{\bf C}_{\bf i}}$ satisfying the following
conditions: \par $(1)$ $R'$ is an algebra over the real field $\bf
R$ under point-wise multiplication and addition;
\par $(2)$ if $x_1,...,x_n\in R$ and $f_1,....,,f_n\in F$, then
$(x_1f_1+...+x_nf_n)\in R'$; \par $(3)$ $R'$ is a normed algebra so
that $ \| x^f \| = \| x \| |f|$ in the case over ${\cal A}_r$ or $
\| x^f \| \le \| x \| |f|$ over $({\cal A}_r)_{{\bf C}_{\bf i}}$ for
each $x\in R$ and $f\in F$, where $x^f := xf$, $~ |f|$ denotes a
norm of $f$ in $F$;
\par $(4)$ the $\bf R$-linear combinations of Form $(2)$ are dense in $R'$;
\par $(5)$ If $x=x_1^{f_1}+...+x_n^{f_n}$ and $\phi $ is a
continuous homomorphism as above, then $ \| x_1\phi (f_1)
+...+x_n\phi (f_n) \| \le \| x \| $.
\par Each homomorphism $\phi $ of $F$ induces $\hat{\phi }(x)$ with the
property $\hat{\phi }(x^f)= x \phi (f)$ and $\hat{\phi }$ will be
called a generated homomorphism.

\par {\bf 72. Theorem.} {\it Let suppositions of \S 71 be satisfied.
Then and element $x\in R'$ has a left inverse if for each generated
homomorphism $\hat{\phi }$ the corresponding element $\hat{\phi
}(x)$ of $R$ has a left inverse in $R$. }
\par {\bf Proof.} Since a homomorphism $\phi $ is $\bf
R$-homogeneous and additive, then it is $\bf R$-linear. Take an
arbitrary maximal ideal $\cal I$ in $R'$. It has the decomposition
\par $(1)$ ${\cal I} = \bigoplus_{j=0}^{2^r-1} {\cal I}_ji_j$, \\ where ${\cal
I}_j$ is either a real or complex algebra isomorphic with ${\cal
I}_k$ for each $0\le j, k \le 2^r-1$. Each $x\in R'$ has the
corresponding element $\theta (x)$ of $H$ (applying Lemma 68 to $R'$
here instead of $R$ in \S 68). \par The algebra $R'$ has the
decomposition $R'={R'}_0i_0\oplus ... \oplus {R'}_mi_m$ induced by
that of $R$ with pairwise isomorphic commutative algebras ${R'}_j$
and ${R'}_k$ either over $\bf R$ or ${\bf C}_{\bf i}$ respectively
for each $k, j$, $~m=2^r-1$. Thus any two elements $a, b\in R'$
quasi-commute and $a=a_0i_0+...+a_mi_m$ and $b=b_0i_0+....+b_mi_m$
with $a_j, b_j\in {R'}_j$ for each $j$. Particularly, elements
$I^f=If$ of $R'$ quasi-commute with each $x^g\in R'$ and hence with
each $b\in R'$. In view of Theorem I.2.81 and Corollary I.2.84
\cite{ludopalglamb} and Lemmas 67 and 70 above the mapping $\theta
(I^f)$ is the continuous algebraic homomorphism from $F$ into ${\cal
A}_r$ or $({\cal A}_r)_{{\bf C}_{\bf i}}$ correspondingly. There
exists a homomorphism $\phi $ so that $\theta (I^f)=J\phi (f)$,
where $J:=\theta (I1)$ is a unit of $H$. Therefore, $\hat{\phi }
(x_1^{f_1}+...+x_n^{f_n})=\theta (x_11)\theta (If_1)+...+\theta
(x_n1)\theta (If_n)=\theta (x_11)\phi (f_1)+...+\theta (x_n1)\phi
(f_n)= \theta ([x_1\phi (f_1)+...+x_n\phi (f_n)]1)$, consequently,
$\theta (x)=\theta (\hat{\phi }(x)1)$ for each
$x=(x_1f_1+...+x_nf_n)\in R'$. From Condition 71$(5)$ and Lemma 70
it follows that $\theta (x)=\theta (\hat{\phi }(x)1)$ for each $x\in
R'$.
\par If $\hat{\phi }(x(t))$ has a left inverse $y$ in $R$, then
$\theta (y1)\theta (x)=\theta (yx)$, consequently, $\theta (yx) =
\theta (\hat{\phi }(yx)1)=\theta (y\hat{\phi }(x)1)=\theta (I1)=J$.
This means that $\theta (x)$ has a left inverse, hence by Lemma 68
$x$ has a left inverse in $R'$.
\par {\bf 73. Corollary.} {\it If suppositions of \S 71 are fulfilled
and for each $\phi (f)$ with $f\in F$ a point $t_0$ exists so that
$\phi (f)=f(t_0)$, then Condition 71$(5)$ can be replaced by $ \|
x(t_0)\| \le \| x \| $ for each $t_0$. Moreover, in the latter
situation an element $x\in R'$ has a left inverse in $R'$, if $x(t)$
has a left inverse in either $R$ or $R_{{\bf C}_{\bf i}}$
correspondingly for each $t$. }
\par {\bf 74. Remark.} If an algebra $F$ has not a unit, then one
can formally adjunct a unit $1$ and consider an algebra ${\bar F} :=
\{ c1+f: ~ c\in Q, ~ f\in F \} $, where either $Q={\cal A}_r$ or
$Q=({\cal A}_r)_{{\bf C}_{\bf i}}$ correspondingly, putting
$|c1+f|^2=|c|^2+|f|^2$ and ${\bar R}' := \{ z=cI1+x: ~ x\in R', ~ c
\in Q \} $ with $ \| z \| ^2 = |c|^2+\|x\|^2$. This standard
construction induces an extended homomorphism either ${\bar {\phi }}
(c1+f)=c+\phi (f)$ or an exceptional homomorphism ${\bar {\phi }}
(c1+f)=c$. If $F$ has not a unit, then statements above can be
applied to $\bar F$ and ${\bar R}'$ so that an element ${\bar x}
=cI1+x$ with $c\ne 0$ may have a left inverse of the form $(bI1+y)$.

\par {\bf 75. Corollary.} {\it Suppose that $\Gamma $ is an additive discrete group so that
$\Gamma = \Gamma_0i_0\oplus ...\oplus \Gamma _mi_m$ with pairwise
isomorphic commutative groups $\Gamma _j$ and $\Gamma _k$ for each
$0\le j, k \le m$ with $m=2^r-1$, $~0\le r \le 3$, while
$G=G_0i_0\oplus ... \oplus G_mi_m$ is an additive group so that
$G_j$ is dual to $\Gamma _j$ with continuous characters $\chi (\beta
,t) = \prod_{k=0}^m \chi _k(\beta _k,t_k)\in S^1$, where $S^1 := \{
u\in {\bf C}_{\bf i}: ~ |u|=1 \} $, $~ \beta =\beta _0i_0+...+\beta
_mi_m$ and $t=t_0i_0+...+t_mi_m$ with $\beta _k\in \Gamma _k$ and
$t_k\in G_k$ for each $k$. Let $$(1)\quad x(t)=\sum_{\beta \in
\Gamma } a_{\beta } \chi (\beta ,t)$$ with $a_{\beta }\in R$ for
each $\beta \in \Gamma $ and
$$(2)\quad \sum_{\beta \in \Gamma } \| a_{\beta } \| <\infty ,$$
then $x\in R'$ has a left inverse in $R'$ if $x(t)$ has a left
inverse in $R_{{\bf C}_{\bf i}}$ for each $t\in G$.}

\par {\bf 76. Corollary.} {\it Let suppositions of Corollary 75 be satisfied,
but with $(2)$ replaced by $$(1)\quad \sum_{\beta \in \Gamma }
e^{q(\beta )} \| a_{\beta } \| <\infty ,$$ where $q(\beta )\in \bf
R$ and \par $(2)$ $q(\alpha + \beta ) \le q(\alpha ) + q(\beta )$
and $q(0)=0$. Then $x\in R'$ is invertible, if $[\sum_{\beta \in
\Gamma } a_{\beta } e^{p(\beta )} \chi (\beta ,t)]$ has a left
inverse in $R_{{\bf C}_{\bf i}}$ for each $t\in G$ with a system of
reals $p(\beta )\in \bf R$ so that \par $(3)$ $p(\alpha + \beta
)=p(\alpha )+p(\beta )$ and $p(0)=0$ and $p(\beta )\le q(\beta )$
for each $\beta \in \Gamma $.}
\par {\bf 77. Corollary.} {\it Let suppositions of Corollary 76 be
satisfied, but let $\Gamma $ be a locally compact group with a
nontrivial nonnegative Haar measure $\lambda $. If $R'$ is formed by
elements of the form:
$$(1)\quad x(t) = \int_{\Gamma } a_{\beta } \chi (\beta ,t)\lambda
(d \beta )$$ with $a_{\beta }\in R$ and
$$(2)\quad \int_{\Gamma } e^{q(\beta )} \| a_{\beta } \|
\lambda (d \beta ) < \infty .$$ If $R'$ has the unit $1(t)=1$ for
each $t\in G$, then $x$ is left invertible, if $x(t)$ has a left
inverse in $R_{{\bf C}_{\bf i}}$ for each $t\in \Gamma $. If $R'$
has not a unit, but $1$ is an adjoint unit as in \S 74, then an
element ${\bar x} =c1+x$ with $c\ne 0$ has a left inverse of the
form $b1+y$, if $[cI+\int_{\Gamma } a_{\beta } e^{p(\beta )} \chi
(\beta ,t) \lambda (d\beta ) ]$ has a left inverse for every $t\in
G$ and each continuous system $p(\beta )$ satisfying Conditions
76$(2,3)$.}
\par {\bf 78. Remark.} Duality theory for locally compact groups is
contained in  \cite{pontr,hew}. Particularly, ${\cal A}_r$ can be
considered as the additive commutative group $({\cal A}_r,+)$. As
the additive group it is isomorphic with ${\bf R}^{2^r}$. The group
of characters of ${\bf R}^n$ is isomorphic with ${\bf R}^n$ for any
natural number $n$ (see \S 23.27(f) in Chapter 6 of the book
\cite{hew}). The Lebesgue measure on the real shadow ${\bf R}^{2^r}$
induces the Lebesgue measure $\lambda $ on ${\cal A}_r$, which is
the Haar measure on $({\cal A}_r,+)$ (see also \S 1).
\par It is possible to consider a dense subgroup $K$ of the total
compact dual group $G$, when $\Gamma $ is discrete. It is sufficient
an existence of a left inverse $y(t)$ of $x(t)$ for each $t\in K$
and that $\sup_{t\in K} \| y(t) \| <\infty $ due to the following
lemma.
\par {\bf 79. Lemma.} {\it Let $x_n$ tend to $x$ in $R$,
when a natural number $n$ tends to the infinity, let also $y_n$ be a
left inverse of $x_n$ for each $n$ and $\sup_n \| y_n \| <\infty $,
then $x\in R$ possesses a left inverse. }
\par {\bf Proof.} From the equality $I-y_nx=I-y_nx_n+y_nx_n-y_nx$ it follows that
$ \| I- y_nx \| \le \| y_n \| \| x_n - x \| $. Then Lemma 53 implies
that a natural number $k$ exists so that $y_nx$ has a left inverse
$z_n$ for each $n\ge k$, consequently, $z_ny_n$ is a left inverse of
$x$ due to the alternativity of the algebra $R$ or using Moufang's
identities.
\par {\bf 80. Corollaries.} {\it Suppose that an algebra $~R$ is over the Cayley-Dickson
algebra ${\cal A}_v$ (see \S 52) and $0\le r\le v$ and $2\le v \le
3$, $~\Gamma = ({\cal A}_r,+)$ (see \S \S 75-78). \par $(1)$. If
$$x(t) =\sum_n a_n e^{(\beta (n), t){\bf i}}\in R'$$ with $a_n\in R$
and $\sum_n \| a_n \| <\infty $, $ ~ (\beta ,t)=Re (\beta t^*) =
\beta _0t_0+...+\beta _mt_m$, where $~ \beta =\beta _0i_0+...+\beta
_mi_m\in ({\cal A}_r,+)$, $~t=t_0i_0+...+t_mi_m\in G = ({\cal
A}_r,+)$, $~m=2^r-1$, moreover, a left inverse $z(t)$ exists for
each $t$ and $ \sup_t \| z(t) \| =C<\infty $, then a function
$$y(t)= \sum_n b_n e^{(\tau (n), t){\bf i}}\in R'$$ exists with
$\sum_n \| b_n \| <\infty $ such that $y(t)x(t)=I$ for each $t$.
\par $(2)$. If $a(\beta )\in R$ is a strongly integrable function
with
$$\int_{{\cal A}_r} \| a(\beta ) \| \lambda (d\beta )<\infty $$ and
if for a nonzero complex number $c\in {\bf C}_{\bf i}\setminus \{ 0
\} $ a function $$[cI + \int_{{\cal A}_r} a(\beta ) e^{(\beta
,t){\bf i}} \lambda (d\beta )]$$ has a left inverse for all $t$ (see
$\lambda $ in \S 78), then a left inverse of the form
$$[qI+\int_{{\cal A}_r} b(\beta ) e^{(\beta ,t){\bf i}} \lambda
(d\beta )]$$ exists with
$$\int_{{\cal A}_r} \| b(\beta ) \| \lambda (d\beta )<\infty .$$}

\par {\bf 81. Corollary.} {\it The algebra $L_q^{n,per}(l_{\infty }
({\bf Z},Y))$ is the saturated subalgebra in the algebra
$L_q^c(l_{\infty } ({\bf Z},Y))$, where $Y$ is a Banach space over
either the quaternion skew field or the octonion algebra ${\cal
A}_v$, $~2\le v\le 3$.}

\par {\bf 82. Theorem.} {\it Suppose that $B\in L_q^{n,per}(l_{\infty }
({\bf Z},Y))$, where $Y$ is a Banach space over either the
quaternion skew field or the octonion algebra ${\cal A}_v$, $~2\le
v\le 3$. Then the following conditions are equivalent:
\par $(1)$ an operator $B$ is invertible;
\par $(2)$ a Fourier transform operator $\hat{B}(M)$ is invertible
for each $M\in {\bf S}^1$ and $x\in Y$.}
\par {\bf Proof.} A real linear Banach subspace $X_k$ is considered, which is linearly
isometrically isomorphic with $l_{\infty } ({\bf Z},Y_k)$ for each
$k\ge 0$. On the other hand, the real span $span_{\bf R} \{ x\in
X_ki_k: ~ k\ge 0 \} $ is dense in $X$. In view of Theorem 51 an
operator $B$ is uniformly $c$-continuous, $B\in L_q^{uc}(l_{\infty }
({\bf Z},Y))$. Therefore, its invertibility on $span_{\bf R} \{ x\in
X_ki_k: ~ k\ge 0 \} $ is equivalent to that of on $X$.
\par Let $B\in L_q^{n,per}(l_{\infty }
({\bf Z},Y))$ be an invertible operator with $D=B^{-1}$. In view of
Corollary 81 the inverse operator $D$ is $n$-periodic as well and
has an absolutely converging Fourier series. From Theorem 51 it
follows that $D$ is uniformly $c$-continuous. Applying Corollary 47
and Proposition 45 we deduce that the Fourier transform operator
$\hat{B}$ is invertible for each $M\in {\bf S}^1$ and $x\in Y$, that
is $\hat{B}(M)\hat{D}(M)x=\widehat{BD}(M)x=I_Yx=x$. Thus
$(1)\Rightarrow (2)$.
\par Vise versa suppose that Condition $(2)$ is fulfilled.
Then by Corollary 81 the mapping $\psi: ~ M\mapsto
(\hat{B}(M))^{-1}$ has an absolutely converging Fourier series: \par
$$(1)\quad \psi (M) = \sum_{k= -\infty }^{\infty } M^k \mbox{ }_kD,
\mbox{ where}$$
$$(2)\quad  \mbox{}_kD := \frac{1}{2\pi } \int_0^{2\pi } e^{-{\bf i}tk}
\psi (e^{{\bf i}t})dt\in L_q(Y\oplus Y{\bf i}) .$$  Put $D= \sum_{k=
-\infty }^{\infty } \mbox{ }_k\bar{D} S(k)$, where as usually $S(k)$
denotes the coordinatewise shift operator on $k$, $S(k)x(l)=x(l+k)$,
while $\mbox{ }_k\bar{D}$ denotes an $1$-ribbon operator, matrix
elements of which on the $k$-th diagonal are equal to $\mbox{}_kD$.
Therefore, the operator $D$ is bounded with $ \| D \| \le \sum_{k=
-\infty }^{\infty } \| \mbox{}_kD \| = c <\infty $. In accordance
with Corollary 47 the operator $D$ is inverse of $B$, i.e.
$D=B^{-1}$.

\par {\bf 83. Corollary.} {\it Let $B\in L_q^{n,per}(l_{\infty }
({\bf Z},Y))$, where $Y$ is a Banach space over either the
quaternion skew field or the octonion algebra ${\cal A}_v$, $~2\le
v\le 3$. Then spectral sets of $\bf B$ and ${\bf B}(D({\cal M}))$
are related by the formula:
\par $\sigma ({\bf B}) = \bigcup_{M\in {\bf S}^1}
\sigma ({\bf B}(D({\cal M}))$, where  ${\bf B}$ denotes the natural
extension of $B$ from $l_{\infty } ({\bf Z},Y)$ onto $l_{\infty }
({\bf Z},Y\oplus Y{\bf i})$.}
\par {\bf Proof.}  The spectral set $\sigma ({\bf B})$ is the complement of
the resolvent set (see Definition I.2.6 \cite{ludopalglamb}), where
$Y\oplus Y{\bf i}$ and $l_{\infty } ({\bf Z},Y\oplus Y{\bf i})$ have
structures of ${\cal A}_v$ Banach spaces as well. In view of
Proposition 45 and Theorem 82 one gets this corollary.

\par {\bf 84. Theorem.} {\it Let a kernel $K$ of a periodic
operator $B$ from \S 1 satisfy the condition:
$$(1)\quad \sup_{t, s}  \| K(t,s) \| =
c_1 <\infty ,$$ where $2\le v \le 3$. Then an operator $A=I-B$ is
invertible if and only if a Fourier transform operator ${\hat A}(M)$
is invertible on $Y\oplus Y{\bf i}$ for each $M\in {\bf S}^1$.}
\par {\bf Proof.} Condition $(1)$ implies that an operator $B$ is bounded
and the integral of \S 1 exists due to the theorem in \S 1 (see also
Chapter 10 in \cite{ferpelruizb}). Take an operator $A=I-B\in
L_q(L^p({\cal A}_w,Y))$, where $p\in [1,\infty ]$. Choose a domain
$V$ in the Cayley-Dickson algebra ${\cal A}_w$ so that $V= \{ z: ~
z\in {\cal A}_w; ~ \forall j ~ z_j\in [0,\omega _j]; ~
z=\sum_{j=0}^{2^w-1} z_ji_j \} $. Then we define an operator $U:
L^p({\cal A}_w,Y)\to l_p({\bf Z}^{2^w},L^p(V,Y))$ by the formula:
$((Ux)(t))({\bar m}) := y_{\bar m}(\tau )$, where $x\in L^p({\cal
A}_w,Y)$ and $y\in l_p({\bf Z}^{2^w},L^p(V,Y))$ are related by the
equation: $y_{{\bar m}}(\tau )=x(t - \sum_j m_j \omega _j i_j)$ with
$\tau _j= t_j - m_j\omega _j\in [0,\omega _j]$ for every $0\le j\le
2^w-1$, where ${\bar m}=(m_0,...,m^{2^w-1})\in {\bf Z}^{2^m}$. This
definition implies that such operator $U$ is an invertible isometry.
\par There exists an operator $Q=UAU^{-1}$, hence $Q\in L_q(l_p({\bf
Z}^{2^w},L^p(V,Y)))$. Evidently $Q$ is invertible if and only if $A$
is such. If $S(\omega )x(t) = x(t+\omega )$ and ${\hat S}({\bar m})
y_{\bar n} = y_{{\bar m} + {\bar n}}$ are shift operators, they
satisfy the equation ${\hat S}({\bar m}) = U S(\sum_j m_j \omega _j
i_j)U^{-1}$. Therefore, the operator $Q$ commutes with each shift
operator ${\hat S}({\bar m})$, where ${\bar m}\in {\bf Z}^{2^w}$.
Thus, this operator $Q$ is $1$-periodic by each $m_j\in {\bf Z}$
(see Definition 33). \par For a function $x\in L^p(V,Y)$ put ${\bar
x}_{\bar m}({\bar k}) := [\prod_{j=0}^{2^w-1} \delta _{m_j,k_j}] x$,
hence ${\bar x}_{\bar m} \in l_p({\bf Z}^{2^w},L^p(V,Y))$. By each
variable $m_j$ a matrix of the operator $Q$ takes the form:
$Q_{k_j,s_j}x = (Q {\bar x}_{\bar s})({\bar k}) = (UAU^{-1}{\bar
x}_{\bar s})({\bar k}) = (UA y_{\bar s})({\bar k})$, when $k_l=s_l$
for each $l\ne j$, where $y_{\bar s}\in l_p({\bf Z}^{2^w},L^p(V,Y))$
is given by the formula: $y_{\bar s}(\tau )=0$ if there exists $j$
so that $\tau _j\notin [s_j\omega _j,(s_j+1)\omega _j)$, while
$y_{\bar s}(\tau )=x(t - \sum_j s_j\omega _ji_j)$ when $\tau _j\in
[s_j\omega _j,(s_j+1)\omega _j)$ for each $j$, where $\tau = t -
\sum_j s_j\omega _ji_j$. This means that a tensor operator $Q^{\bar
s}_{\bar k}$ is defined: $Q^{\bar s}_{\bar k}= (Q{\bar x}_{\bar
s})({\bar k})$. Then the function $(UA y_{\bar s})({\bar k})\in
L^p(V,Y)$ takes the form:
$$(2)\quad (UA y_{\bar s})({\bar k})(\tau ) = x(\tau )$$
$$ - \mbox{}_{\sigma }
\int_{\gamma ^{\alpha }(b)|_{b_0\in [s_0\omega _0,(s_0+1)\omega
_0]}}...\mbox{}_{\sigma }\int_{\gamma ^{\alpha }(b)|_{b_u\in
[s_u\omega _u,(s_u+1)\omega _u]}} K((\tau + \sum_j k_j\omega
_ji_j),b)$$
$$x(b - \sum_j s_j\omega _ji_j)) db_0...db_u$$  $$ =
x(\tau ) - \mbox{}_{\sigma } \int_{\gamma ^{\alpha }(b)|_{b_0\in
[s_0\omega _0,(s_0+1)\omega _0]}}...\mbox{}_{\sigma }\int_{\gamma
^{\alpha }|_{b_u\in [s_u\omega _u,(s_u+1)\omega _u]}} K((\tau +
\sum_j k_j\omega _ji_j),$$  $$(b+\sum_j s_j\omega _ji_j)) x(b)
db_0...db_u$$
$$=x(\tau ) - \mbox{}_{\sigma } \int_{\gamma ^{\alpha }(b)|_{b_0\in
[s_0\omega _0,(s_0+1)\omega _0]}}...\mbox{}_{\sigma }\int_{\gamma
^{\alpha }(b)|_{b_u\in s_u\omega _u,(s_u+1)\omega _u]}} K(\tau ,$$
$$ b +\sum_j (s_j-k_j)\omega _ji_j) x(b) db_0...db_u$$
in accordance with Conditions 1$(5,6)$, where $u := 2^w-1$. This
implies that a tensor operator takes the form:
$$(3)\quad Q^{\bar s}_{\bar k}= x(\tau ) - \mbox{}_{\sigma }
\int_{\gamma ^{\alpha }(b)|_{b_0\in [s_0\omega _0,(s_0+1)\omega
_0]}}...\mbox{}_{\sigma }\int_{\gamma ^{\alpha }(b)|_{b_u\in
[s_u\omega _u,(s_u+1)\omega _u]}} K(\tau ,$$  $$ b +\sum_j
(s_j-k_j)\omega _ji_j) x(b) db_0...db_u$$ for each ${\bar s}, {\bar
k}\in {\bf Z}^{2^w}$. Elements of this  tensor depend only on the
difference ${\bar s}-{\bar k}$, so it is possible to put $Q^{\bar
s}_{\bar k}=Q_{{\bar s}-{\bar k}}$.
\par The operators $U$, $A$ and $U^{-1}$ are $c$-continuous by each $s_j$,
consequently, the operator $Q_{\bar s}$ is also $c$-continuous by
each $s_j$, where $j = 0, 1, ..., 2^w-1$. Applying Theorem 82 we
obtain the statement of this theorem.

\par {\bf 85. Corollary.} {\it Let an operator $B$ satisfy
conditions of Theorem 84, then a spectral set is $\sigma ({\bf Q}) =
\bigcup_{M\in {\bf S}^1} \sigma ({\bf Q}(D({\cal M})))$, where ${\bf
Q}$ denotes the natural extension of $Q$ from $l_{\infty } ({\bf
Z}^{2^w},L^{\infty }(V,Y))$ onto $l_{\infty } ({\bf
Z}^{2^w},L^{\infty }(V,Y\oplus Y{\bf i}))$.}
\par {\bf Proof.} This follows from Theorem 84 and Corollary 83 applying
the Fourier transform by each variable, since ${\bf Q}\in
L_q^{1,per}(l_{\infty } ({\bf Z}^{2^w},L^{\infty }(V,Y)\oplus
L^{\infty }(V,Y){\bf i}))$ due to Condition 84$(1)$ and the latter
Banach space over the Cayley-Dickson algebra ${\cal A}_v$ is
isomorphic with $L_q^{1,per}(l_{\infty } ({\bf Z}^{2^w},L^{\infty
}(V,Y\oplus Y{\bf i})))$.


\begin{thebibliography}{399}

\bibitem{ablsigb} M. J. Ablowitz, H. Segur. "Solitons and the inverse
scattering transform" (SIAM: Philadeplhia, 1981).

\bibitem{annmbochph42} S. Bochner, R.S. Phillips. "Absolutely
convergent Fourier expansions for non-commutative normed rings".
Annals of Mathem. {\bf 43: 3} (1942), 409-418.

\bibitem{baez} J.C. Baez. "The octonions". Bull. Amer.
Mathem. Soc. {\bf 39: 2} (2002), 145-205.

\bibitem{brdeso} F. Brackx, R. Delanghe, F. Sommen.
"Clifford analysis" (London: Pitman, 1982).

\bibitem{dickson} L.E. Dickson. "The collected mathematical papers".
Volumes 1-5 (Chelsea Publishing Co.: New York, 1975).

\bibitem{emch} G. Emch. "M$\grave e$chanique quantique quaternionienne et
Relativit$\grave e$ restreinte". Helv. Phys. Acta {\bf 36} (1963),
739-788.

\bibitem{ferpelruizb} J.C. Ferrando, M. L\'opez Pellicer, L.M.
S\'anchez Ruiz. "Metrizable barreled spaces" (Longman Group Ltd:
Harlow, 1995).

\bibitem{fihteng} G.M. Fihtengolz. "Course of differential and
integral calculus", 8-th Edition, V. 1-3 (Moscow: Fizmatlit, 2003).

\bibitem{gilmurr} J.E. Gilbert, M.A.M. Murray.
"Clifford algebras and Dirac operators in harmonic analysis". Cambr.
studies in advanced Mathem. {\bf 26} (Cambr. Univ. Press: Cambridge,
1991).

\bibitem{girard} P.R. Girard. "Quaternions, Clifford algebras and
relativistic Physics" (Birkh\"auser: Basel, 2007).

\bibitem{guesprqa} K. G\"urlebeck, W. Spr\"ossig. "Quaternionic
analysis and elliptic boundary value problem" (Birkh\"auser: Basel,
1990).

\bibitem{guetze} F. G\"ursey, C.-H. Tze. "On the role of
division, Jordan and related algebras in particle physics" (World
Scientific Publ. Co.: Singapore, 1996).

\bibitem{hew} E. Hewitt, K.A. Ross. "Abstract harmonic analysis"
(Berlin: Springer, 1979).

\bibitem{kansol} I.L. Kantor, A.S. Solodovnikov.
"Hypercomplex numbers" ( Springer-Verlag: Berlin, 1989).

\bibitem{krausryan} R.S. Krausshar, J. Ryan. "Some conformally
flat spin manifolds, Dirac operators and automorphic forms". J.
Math. Anal. Appl. {\bf 325} (2007), 359-376.

\bibitem{kravchot} V.V. Kravchenko. "On a new approach for solving
Dirac equations with some potentials and Maxwell's sytem in
inhomogeoneous media". Operator Theory {\bf 121} (2001), 278-306.

\bibitem{kuznpr} V.V. Kuznetzov. "Spectral properties of periodic
integral operators". Prepr. {\bf 11} (2000), 1-32 (RAN DVO:
Vladivostok, 2000).

\bibitem{ludoyst} S.V. Ludkovsky, F. van Oystaeyen.
"Differentiable functions of quaternion variables". Bull. Sci. Math.
(Paris). Ser. 2. {\bf 127} (2003), 755-796.

\bibitem{ludfov} S.V. Ludkovsky. "Differentiable functions of
Cayley-Dickson numbers and line integration". J. of Mathem. Sciences
{\bf 141: 3} (2007), 1231-1298.

\bibitem{lujmsalop} S.V. Ludkovsky.  "Algebras of operators in Banach
spaces over the quaternion skew field and the octonion algebra". J.
Mathem. Sciences {\bf 144: 4} (2008), 4301-4366.

\bibitem{lufjmsrf} S.V. Ludkovsky. "Residues of functions of octonion
variables". Far East Journal of Mathematical Sciences (FJMS), {\bf
39: 1} (2010), 65-104.

\bibitem{ludancdnb} S.V. Ludkovsky. "Analysis over Cayley-Dickson
numbers and its applications" (LAP Lambert Academic Publishing:
Saarbr\"ucken, 2010).

\bibitem{luspraaca} S.V. Ludkovsky, W. Sproessig.  "Ordered
representations of normal and super-differential operators in
quaternion and octonion Hilbert spaces". Adv. Appl. Clifford Alg.
{\bf 20: 2} (2010), 321-342.

\bibitem{ludspr} S.V. Ludkovsky, W. Spr\"ossig.
"Spectral theory of super-differential operators of quaternion and
octonion variables", Adv. Appl. Clifford Alg. {\bf 21: 1} (2011),
165-191.

\bibitem{ludspr2} S.V. Ludkovsky, W. Spr\"ossig.
"Spectral representations of operators in Hilbert spaces over
quaternions and octonions", Complex Variables and Elliptic
Equations, online, DOI:10.1080/17476933.2010.538845, 24 pages
(2011).

\bibitem{ludvhpde} S.V. Ludkovsky.  "Integration of vector
hydrodynamical partial differential equations over octonions".
Complex Variables and Elliptic Equations, online,
DOI:10.1080/17476933.2011.598930, 31 pages (2011).

\bibitem{ludcmft12} S.V. Ludkovsky. "Line integration of Dirac operators
over octonions and Cayley-Dickson algebras". Computational Methods
and Function Theory, {\bf 12: 1} (2012), 279-306.

\bibitem{ludopalglamb} S.V. Ludkovsky. "Operator algebras over Cayley-Dickson numbers"
(LAP LAMBERT Academic Publishing AG $\&$ Co. KG: Saarbr\"ucken,
2011).

\bibitem{ludunbnormopla12} S.V. Ludkovsky. "Unbounded normal
operators in octonion Hilbert spaces and their spectra", Los Alamos
Nat. Lab., math.FA/1204.1554 (2012), 49 pages.

\bibitem{oystaey} F. van Oystaeyen. "Algebraic geometry for
associative algebras". Series "Lect. Notes in Pure and Appl.
Mathem." {\bf 232} (Marcel Dekker: New York, 2000).

\bibitem{pontr} L.S. Pontrjagin. "Continuous groups"
(Moscow: Nauka, 1984).

\bibitem{schafb} R.D. Schafer. "An introduction to non-associative
algebras" (Academic Press: New York, 1966).

\bibitem{wienannm32} N. Wiener. "Tauberian theorems". Abnnals of
Mathematics. {\bf 33: 1 } (1932), 1-100.

\bibitem{zelditch09} S. Zelditch. "Inverse spectral problem for
analytic domains, II: ${\bf Z}_2$-symmetric domains". Advances in
Mathematics {\bf 170: 1} (2009), 205-269.

\end{thebibliography}
\end{document}